\documentclass[amsmath,10pt,a4paper]{article}
\usepackage{amsfonts, amssymb, amsmath, amscd}
\usepackage[french]{babel}
\usepackage{amssymb}
\newtheorem{thm}{Th\'eor\`eme}[section]

\newtheorem{lem}[thm]{Lemme}
\newtheorem{prop}[thm]{Proposition}
\newtheorem{cor}[thm]{Corollaire}

\numberwithin{equation}{section}

\title{Produits dans la cohomologie des vari\'et\'es arithm\'etiques~: quelques calculs sur les s\'eries th\^eta.}

\author{N. Bergeron}

\date{}

\begin{document}

\maketitle

\begin{footnotesize}
\begin{center}
{\bf Abstract}
\end{center}
For abelian varieties $A$, in the most interesting cohomology theories $H^* (A)$ is the exterior algebra of $H^1(A)$. In this paper we study a weak generalization
of this in the case of arithmetic manifolds associated to orthogonal or unitary groups. In this latter case recall that arithmetic manifolds associated to 
standard unitary groups $U(p,q)$ ($p\geq q$) over a totally real numberfield have vanishing cohomology in degree $i=1, \ldots , q-1$ and that,  
following earlier works of Kazhdan and Shimura, Borel and Wallach constructed in \cite{BorelWallach} non zero degree $q$ cohomology classes.
These cohomology classes arise as theta series.  After generalizing the construction of these theta series. 
We prove that arbitrary (up to the obvious obstructions) cup-products of these theta series and their complex conjugates virtually non vanish, i.e.
``up to Hecke translate'', in the cohomology ring. This fits inside the, partly conjectural, picture drawn in \cite{Lefschetz}. 
\end{footnotesize}

\section{Introduction}

Cet article fait suite \`a \cite{Lefschetz} dans lequel nous formulons, et d\'emontrons dans un grand nombre de cas, des conjectures concernant la cohomologie des
vari\'et\'es arithm\'etiques associ\'ees aux groupes orthogonaux et unitaires. De la m\^eme mani\`ere nous ne consid\'erons ici que le cas d'un espace sym\'etrique
$D^+$ associ\'e au groupe $O(p,q)$ ou $U(p,q)$. Soit plus pr\'ecisemment $k$ un corps de nombres totalement r\'eel (resp. une extension quadratique imaginaire
d'un corps de nombres totalement r\'eel). Notons $\sigma_1 , \ldots , \sigma_{\mu}$ 
les diff\'erents plongements archim\'ediens (resp. plongements archim\'ediens modulo conjugaison) de $k$.
Dans le second cas, notons $k_0$ le sous-corps de nombres totalement r\'eel fix\'e par la conjugaison complexe. Soit $V_k$ un espace vectoriel sur $k$ et 
$(,)$ une forme quadratique (resp. hermitienne) non d\'eg\'en\'er\'ee et anisotrope 
sur $V_k$ de signature $(p,q)$ en une place archim\'edienne de $k$ et d\'efinie positive en les 
autres places. Soit $V$ l'espace vectoriel des points r\'eels de l'espace vectoriel obtenu \`a partir de $V_k$ par restriction des scalaires de $k$ \`a ${\Bbb Q}$ (resp.
de $k_0$ \`a ${\Bbb Q}$). On a alors un isomorphisme d'espaces vectoriels r\'eels
$$V \cong \oplus_{j=1}^{\mu} V^{(j)},$$
o\`u $V^{(j)}$ est le compl\'et\'e de $V_K$ relativement au plongement $\sigma_j$. Nous supposerons que la forme $(,)_1$ induite par $(,)$ sur $V^{(1)}$ est de signature 
$(p,q)$ et donc que les formes $(,)_j$ induites par $(,)$ sur les $V^{(j)}$, pour $j\geq 2$, sont toutes d\'efinies positives. Par abus de notation nous noterons
\'egalement $(,)$ la forme induite sur $V$. Cette forme est somme directe orthogonale des formes $(,)_j$. Notons finalement $G^{(j)}$ le sous-groupe de ${\rm Aut} (V^{(j)})$
constitu\'e des isom\'etries de $(,)_j$, et $G= \prod_j G^{(j)}$. Le groupe $G$ est isomorphe au groupe des points r\'eels du groupe r\'eductif sur ${\Bbb Q}$ obtenu,
par restriction des scalaires de $k$ \`a ${\Bbb Q}$ (resp. $k_0$ \`a ${\Bbb Q}$), \`a partir du groupe des isom\'etries de la forme $(,)$ sur $V_k$. Remarquons
que $G^{(1)} \cong O(p,q)$ (resp. $U(p,q)$) et $G^{(j)} \cong O(m)$ (resp. $U(m)$) pour $j\geq 2$.

Le groupe $G$ est d\'efini sur ${\Bbb Q}$. Un {\it sous-groupe de congruence} de $G({\Bbb Q})$ est un sous-groupe de la forme $\Gamma = G({\Bbb Q}) \cap K$,
o\`u $K$ est un sous-groupe compact ouvert du groupe $G({\Bbb A}_f)$ des points ad\`eliques finis de $G$.

Dans cet article on \'etudie les quotients (compacts, puisque $G$ est anisotrope) $S(\Gamma )= \Gamma \backslash D^+$,
o\`u $\Gamma \subset G({\Bbb Q})$ est un sous-groupe de congruence; ces quotients s'identifient aux composantes
connexes de $G({\Bbb Q}) \backslash G({\Bbb A}) / K_{\infty} K = G({\Bbb Q}) \backslash (D \times G({\Bbb A}_f ) )/K$,
o\`u $K \subset G({\Bbb A}_f )$ est un sous-groupe compact ouvert, $K_{\infty}$ est la pr\'eimage d'un sous-groupe compact maximal de la composante connexe de l'identit\'e 
$G^{{\rm ad}}({\Bbb R})^+$ du groupe adjoint $G^{{\rm ad}} ({\Bbb R})$ de $G({\Bbb R})$ et 
$D = G({\Bbb R}) / K_{\infty}$ est une r\'eunion (finie, disjointe) d'espaces sym\'etriques isom\'etriques \`a $D^+$.
Plus pr\'ecis\'ement, d\'esignons par $G_f$ l'adh\'erence de $G({\Bbb Q})$ dans le groupe $G({\Bbb A}_f )$. Un sous-groupe 
de congruence $\Gamma \subset G({\Bbb Q})$ s'\'ecrit $\Gamma = G({\Bbb Q}) \cap K$ o\`u $K$ est l'adh\'erence de
$\Gamma$ dans $G_f$ et,
\begin{eqnarray} \label{quotient}
S(\Gamma ) = \Gamma \backslash D^+ = G({\Bbb Q}) \backslash (D^+ \times G_f ) / K = S(K).
\end{eqnarray}
On s'interesse ici \`a la cohomologie
(\`a coefficients complexes) $H^* (S(\Gamma) )$ de ces quotients.

Une fois donn\'e deux sous-groupes de congruence $\Gamma ' \subset \Gamma \subset G({\Bbb Q})$, on obtient un
rev\^etement fini
$$S(\Gamma ' ) \rightarrow S(\Gamma )$$
qui induit un morphisme injectif
$$H^* (S(\Gamma ) ) \rightarrow H^* (S(\Gamma ' ) )$$
en cohomologie. Les groupes de cohomologies $H^* (S(\Gamma ) )$ forment donc un syst\`eme inductif index\'e par les
sous-groupes de congruence $\Gamma \subset G({\Bbb Q})$.  En passant \`a la limite (inductive) on d\'efinit
\begin{eqnarray} \label{HSh}
H^* (Sh^0 G) = \lim_{\rightarrow_\Gamma} H^* (S(\Gamma)).
\end{eqnarray}
La notation ci-dessus provient de ce que lorsque l'espace $D^+$ est hermitien, on appelle {\it vari\'et\'e de Shimura connexe}
l'espace topologique
\begin{eqnarray} \label{Sh}
Sh^0 G = \lim_{\leftarrow_\Gamma} S(\Gamma)  = G({\Bbb Q}) \backslash (D^+ \times G_f ).
\end{eqnarray}
On peut consid\'erer sa cohomologie de C\v{e}ch et il est d\'emontr\'e dans \cite{Rohlfs} que celle-ci
co\"{\i}ncide avec (\ref{HSh}).
Pour ce qui nous concerne, il sera suffisant de consid\'erer que $H^* (Sh^0 G)$ n'est qu'une
notation pour la limite inductive (\ref{HSh}).

\medskip

La cohomologie d'une vari\'et\'e ab\'elienne est une alg\`ebre ext\'erieure sur le $H^1$. Dans le cas de 
l'alg\`ebre $H^* (S(\Gamma ))$ on ne peut esp\'erer un \'enonc\'e aussi fort, on peut n\'eanmoins se demander si \'etant donn\'ees deux classes 
$\alpha$ et $\beta$ dans $H^* (S(\Gamma ))$ il est possible de s'assurer que leur cup-produit est non nul.

Cette question est trop na\"{\i}ve. On s'aper\c{c}oit rapidement que pour qu'il y ait un espoir d'y r\'epondre positivement il faut affaiblir la conclusion.
Les r\'esultats de non-annulation que nous avons en vue ne sont que des r\'esultats
{\it virtuels} au sens suivant~: donn\'ees deux classes $\omega$, $\eta$ de cohomologie sur la vari\'et\'e $S(\Gamma)$,
on ne peut affirmer que $\omega \wedge \eta \neq 0$, mais seulement que $(p^* \omega \wedge C'*p^* \eta ) \neq 0$.
Ici $p :S(\Gamma ') \rightarrow S(\Gamma )$ est un rev\^etement galoisien de $S(\Gamma )$, et $C$ est une correspondance de Hecke sur $S(\Gamma ')$. 
En passant \`a la limite inductive (\ref{HSh}), il s'agit donc de trouver des conditions d'alg\`ebre lin\'eaire sur deux classes $\omega$ et $\eta \in 
H^* (Sh^0 G)$ pour qu'il existe un \'el\'ement $g \in G({\Bbb Q})$ tel que le cup-produit $\omega \wedge g(\eta )$ soit non nul dans $H^* (Sh^0 G)$.

Dans \cite{Lefschetz} nous rassemblons sous l'intitul\'e ``propri\'et\'es de Lefschetz automorphes'' une s\'erie de r\'esultats et conjectures qui concernent ce probl\`eme
(ainsi que ceuxi de l'injectivit\'e {\it virtuelle} 
des applications de restriction ou de rel\`evement dans les vari\'et\'es arithm\'etiques associ\'ees aux 
groupes unitaires et orthogonaux). Le but de cet article est d'obtenir des r\'esultats plus complets 
dans le cas de certaines classes sp\'eciales obtenues comme s\'eries th\^eta. 

\paragraph{Classes de cohomologie obtenues comme s\'eries th\^eta.}
\'Etant donn\'ees deux entiers naturels $r$ et $s$, consid\'erons le groupe
symplectique $Sp(2r)$ (resp. le groupe unitaire $U(r,s)$). La paire de groupes r\'eels $(Sp(2r) , O(p,q))$ (resp. $(U(r,s) , U(p,q))$) est une paire r\'eductive 
duale, au sens de Howe \cite{Howe}, dans le groupe symplectique $Sp(2r(p+q))$ (resp. $Sp(2(r+s)(p+q))$). Notons $\omega$ la repr\'esentation de 
Weil du groupe m\'etaplectique correspondant et ${\cal H}$ son module d'Harish-Chandra associ\'e. Dans un premier temps nous construisons explicitement
des \'el\'ements 
$$\varphi^{(rq,sq)} \in \left[ {\cal H}  \otimes \Omega^{(r+s)q} (D) \right]^{G^{(1)}}$$
o\`u $s=0$ dans le cas orthogonal, $\Omega^* (D)$ est l'alg\`ebre des formes diff\'erentielles lisses sur $D$ et le groupe 
$G^{(1)}$ agit sur ${\cal H}$ via la repr\'esentation
de Weil $\omega$. Dans le cas orthogonal ou lorsque $r=s$ dans le cas unitaire, les formes que nous construisons co\"{\i}ncident avec celles construites
par Kudla et Millson dans \cite{KudlaMillson1}. Notre mani\`ere d'y parvenir est compl\`etement diff\'erente. Le principal r\'esultat relatif \`a cette construction est le 
th\'eor\`eme suivant.

\begin{thm} \label{T1}
Les formes diff\'erentielles $\varphi^{(rq,sq)}$ sont ferm\'ees vues comme formes diff\'erentielles sur $D$, {\it i.e.}
$$d \varphi^{(rq,sq)} = 0  $$ 
o\`u $d$ est la diff\'erentielle du complexe calculant la $(\mathfrak{g} , K)$-cohomologie du $G^{(1)}$-module ${\cal H}$.
Elles sont non nulles pour $0 \leq r \leq p$ (et $0 \leq s \leq p$) et v\'erifient les propri\'et\'es suivantes.
\begin{enumerate}
\item Les formes sont compatibles avec le cup-produit~:
$$\varphi^{(r_1 q, s_1 q)} \wedge \varphi^{(r_2 q , s_2 q)} = \varphi^{((r_1 + r_2) q , (s_1 + s_2) q)}  , $$
o\`u $\varphi^{(rq,sq)} =0$ si $r$ ou $s$ est strictement sup\'erieur \`a $p$.
\item Les formes sont, \`a multiplication par une fonction non nulle pr\`es, 
compatibles avec les restrictions de $O(p,q)$ vers $O(p-l,q)$, $U(p,q)$ vers $U(p-l,q)$ ou encore de $U(p,q)$ vers $O(p,q)$.
\item La composante {\it fortement primitive}, autrement dit la composante de $K$-type 
$$(\underbrace{q, \ldots , q}_{r \; {\rm fois}} , 0 , \ldots , 0 , \underbrace{-q, \ldots , -q}_{s \; {\rm fois}}) \otimes (0, \ldots , 0), $$
de 
$$[\varphi^{(rq,sq)} ] \in H^{(r+s)q} (\mathfrak{g} , K ; {\cal H}) $$
est non nulle lorsque $r$ et $s \leq p/2$.
\end{enumerate}
\end{thm}

Pour un \'enonc\'e plus d\'etaill\'e, nous renvoyons \`a la section 3. Retournons maintenant \`a la situation globale du d\'ebut de l'introduction.  

Consid\'erons $V_k '$ un espace vectoriel de dimension fini sur $k$, $(,)'$ une forme symplectique (resp. anti-hermitienne) non d\'eg\'en\'er\'ee sur $V_k '$
et $G'$ le groupe r\'eductif sur ${\Bbb Q}$ obtenu, par 
restriction des scalaires de $k$ (resp. $k_0$) \`a ${\Bbb Q}$, \`a partir du groupe des isom\'etries de $(.)'$. Supposons $G^{' (1)} \cong Sp(2r)$ (resp.
$G^{ ' (1)} \cong U(r,s)$) et $V_k '$ d'indice de Witt maximal (\'egal \`a $|r-s|$) sous cette condition.

Le $k$-espace vectoriel $W_k = V_k \otimes_k V_k '$ est naturellement muni de la forme symplectique 
$$\langle , \rangle = {\rm tr}_{k/k_0} \left( (,) \otimes \overline{(,)'} \right),$$
o\`u $k_0 = k$ dans le cas orthogonal et $x \mapsto \overline{x}$ d\'esigne l'involution de Galois de $k/k_0$. Soit $\widetilde{G}' ({\Bbb A})$
le rev\^etement non trivial (m\'etaplectique) \`a deux feuillets de $G' ({\Bbb A})$ 
(resp. le rev\^etement non trivial \`a deux feuillets qui est trivial au-dessus du groupe sp\'ecial unitaire), cf. Weil \cite{Weil}.

Un caract\`ere non trivial de ${\Bbb A}_{k_0} /k_0$~\footnote{Ici ${\Bbb A}_{k_0}$ d\'esigne les ad\`eles de $k_0$.} \'etant fix\'e,  il existe une repr\'esentation de Weil
$\omega$ de $G ({\Bbb A}) \cdot \widetilde{G}' ({\Bbb A})$. Une {\it polarisation} $W=X \oplus Y$, o\`u $X$ et $Y$ sont des sous-espaces totalement isotropes maximaux
de $W$, donne lieu \`a une r\'ealisation de $\omega$ dans $L^2 (X)$ connue sous le nom de mod\`ele de Schr\"odinger pour $\omega$, cf. Gelbart \cite{Gelbart}.
L'espace de ses vecteurs lisses est l'espace de Bruhat-Schwartz ${\cal S} (X({\Bbb A})) \cong {\cal S} (X_{\infty} ) \otimes  {\cal S} (X({\Bbb A}_f ))$.
Toute fonction $\varphi \in {\cal S} (X({\Bbb A}_f ))$ peut \^etre compl\`et\'ee en une forme de Schwartz globale $\widetilde{\varphi}$ sur $X({\Bbb A})$,
en prenant $\varphi^{(rq,sq)}$ \`a la premi\`ere place archim\'edienne et en prenant la gaussienne qui repr\'esente le vide dans le mod\`ele de Shr\"odinger 
en toutes les autres places archim\'ediennes. La forme $\widetilde{\varphi} \in \left[ {\cal S} (X({\Bbb A})) \otimes \Omega^{(r+s)q} (D) \right]^{G({\Bbb R})}$ et
si $\varphi$ est $K$-invariante, alors pour $g' \in \widetilde{G} ' ({\Bbb A})$ et $g \in G({\Bbb A}_f )$, la s\'erie th\^eta
\begin{eqnarray} \label{thetaserie}
\theta (g,g', \varphi) = \sum_{x \in X(k_0 )} (\omega (gg') \widetilde{\varphi})(x)  
\end{eqnarray}
d\'efinit une forme ferm\'ee $\theta (g' , \varphi )$ de degr\'e $(r+s)q$ sur 
$$Sh (G)_{K} : = G({\Bbb Q})\backslash (D \times G({\Bbb A}_f ) )/ K.$$ 
En restriction \`a la composante connexe $S(\Gamma ) = \Gamma \backslash D^+$,
celle-ci d\'efinit une classe de cohomologie $[\theta (g' , \varphi )] \in H^{(r+s)q} (Sh^0 G)$. Nous dirons de cette classe de cohomologie qu'elle
est d\'efinie par une s\'erie th\^eta et qu'elle est de bidegr\'e $(rq,sq)$.

Ce sont les s\'eries th\^eta du titre. Dans le cas unitaire et pour $(r,s)=(1,0)$, ces s\'eries sont (essentiellement celles) consid\'er\'ees par Kazhdan \cite{Kazhdan},
Shimura \cite{Shimura} et Borel et Wallach \cite{BorelWallach}. Ils montrent dans ce cas que 
\begin{eqnarray} \label{BW}
[\theta (g' , \varphi ) ] \neq 0 \mbox{ dans } H^q (Sh^0 G).
\end{eqnarray} 
Toujours dans le cas unitaire, Anderson consid\`ere dans \cite{Anderson} de telles s\'eries th\^eta avec $s=0$; il montre \'egalement que 
les classes de cohomologie obtenues sont non nulles. 

\paragraph{\'Enonc\'es des r\'esultats.}
Il d\'ecoule facilement du th\'eor\`eme \ref{T1} et de la construction que l'espace de ces s\'eries th\^eta est stable par cup-produit et restriction.
Le principal th\`eme du texte est alors de ``tester'' les propri\'et\'es de Lefschetz conjectur\'ees dans \cite{Lefschetz} sur les classes de cohomologie d\'efinies par ces
s\'eries th\^eta. Relativement aux propri\'et\'es attendues pour le cup-produit nous montrons le th\'eor\`eme suivant.

Si le groupe $G({\Bbb Q})$ (ou $G_f$) agit naturellement sur $H^*(Sh^0 G)$ et pr\'eserve le sous-espace des classes d\'efinies par des s\'eries th\^eta, remarquons que
l'on a m\^eme une action de tout le groupe $G({\Bbb A}_f)$ sur les classes d\'efinies par des s\'eries th\^eta~: celle induite par
l'action de $G({\Bbb A}_f )$ sur $\varphi \in {\cal S} (X({\Bbb A}_f ))$.

\begin{thm} \label{T2}
Soient $[\theta (g_i ' , \varphi_i )] \in H^* (Sh^0 G)$, $i=1,2$, deux classes de cohomologie d\'efinies par des s\'eries th\^eta et de bidegr\'es respectifs 
$(r_i q, s_i q)$, $i=1,2$. Supposons $r_1+r_2 \leq p$ et $s_1 +s_2 \leq p$, il existe alors un \'el\'ement $g \in G({\Bbb A}_f)$ tel que
$$[\theta (g_1 ' , \varphi_1 )] \wedge g \left( [\theta (g_2 ' , \varphi_2 )] \right) \neq 0 \; \mbox{ dans } H^{(r_1+r_2+s_1+s_2) q} (Sh^0 G).$$
\end{thm}

Il est ainsi plaisant de remarquer que, dans le cas unitaire et partant d'une classe de Borel-Wallach (\ref{BW}) $\eta \in H^{(q,0)} (Sh^0 G)$, 
il est possible en formant des cup-produits $\eta_1 \wedge \ldots \wedge \eta_r \wedge \overline{\eta}_1 \wedge \ldots \wedge \overline{\eta}_s$, o\`u les $\eta_i$ (resp.
$\overline{\eta}_j$) sont des translat\'es de Hecke de $\eta$ (resp. $\overline{\eta}$), d'obtenir des classes non nulles dans $H^{(rq,sq)} (Sh^0 G)$
pour $r,s \leq p$. 

\medskip

On peut passer du groupe unitaire au groupe orthogonal \`a l'aide du th\'eor\`eme suivant. Remarquons que si $k/k_0$ est une extension quadratique totalement
imaginaire, la restriction de $(,)$ au sous-espace $V_{k_0}$ de $V_k$ est une forme quadratique. Il correspond \`a tout ceci un 
plongement de $k_0$-groupes~: $O(V_{k_0}) \subset U(V_k )$. Notons $H$ et $G$ les ${\Bbb Q}$-groupes obtenus par restriction des scalaires de $k_0$ \`a
${\Bbb Q}$ et $D_H^+$ et $D_G^+$ les espaces sym\'etriques associ\'es.

\begin{thm} \label{T3}
Soit $\theta (g' , \varphi )$ une s\'erie th\^eta associ\'ee au groupe unitaire $G$ et de bidegr\'e $(rq,0)$ avec $r\leq p$. Alors, il existe un \'el\'ement $g \in G({\Bbb A}_f)$ tel que
la restriction de la forme ferm\'ee $g \left(  \theta (g' , \varphi ) \right)$ sur $D_G^+$ \`a l'espace sym\'etrique $D_H^+$ associ\'e
au groupe orthogonal $H$ d\'efinisse une classe de cohomologie {\bf non nulle} dans $H^{rq} (Sh^0 H)$ (elle aussi d\'efinie par une s\'erie th\^eta).
\end{thm}

En particulier, partant d'une classe de Borel-Wallach (\ref{BW}) dans $H^{(q,0)} (Sh^0 G)$, il est possible, en formant des cup-produits puis en les restreignant
\`a $Sh^0 H$, d'obtenir des classes non nulles dans $H^{rq} (Sh^0 H)$ pour $r \leq p$.

\paragraph{Organisation de l'article.}
Dans une premi\`ere section on fait quelques rappels concernant la repr\'esentation de l'oscillateur harmonique. Les r\'esultats que nous rappelons 
sont tous classiques et d\^us \`a Weil, Cartier, Kashiwara et Vergne, Howe et Kudla.

La seconde section est consacr\'ee \`a la construction des classes $\varphi^{(rq,sq)}$. Celle-ci repose sur la construction de $\varphi^{(q,0)}$ qui est essentiellement
due \`a Borel et Wallach et que l'on explicite \`a l'aide des travaux de Kashiwara et Vergne rappel\'es dans la section pr\'ec\'edente. La consid\'eration de paires duales
en balance, au sens de Kudla, permet de former des ``cup-produits'' de $\varphi^{(q,0)}$ et $\varphi^{(0,q)} = \overline{\varphi^{(q,0)}}$ et d'obtenir nos formes
$\varphi^{(rq,sq)}$ de mani\`ere explicite dans un mod\`ele de Fock de la repr\'esentation de l'oscillateur harmonique. En passant dans un mod\`ele
de Schr\"odinger appropri\'e on montre comment r\'eobtenir les classes de Schwartz de Kudla et Millson \cite{KudlaMillson1}.

Dans la troisi\`eme et derni\`ere section, on globalise ces constructions pour obtenir les s\'eries th\^eta du titre et on d\'emontre les th\'eor\`emes \ref{T2} et \ref{T3}.
La d\'emonstration repose sur les travaux de Kudla et Millson \cite{KudlaMillson2} et \cite{KudlaMillson3} tels qu'appliqu\'es par Kudla dans \cite{Kudla2} et dans le
cas du groupe $O(2,n)$. Notons que l'on a besoin d'un cas de la formule de Siegel-Weil pour les groupes unitaires r\'ecemment d\'emontr\'e par Ichino \cite{Ichino}.

\section{La repr\'esentation de l'oscillateur harmonique}

\paragraph{1.} Soit $(W, B )$ un espace symplectique non d\'eg\'en\'er\'e de dimension 
$2N$ sur ${\Bbb R}$. Le groupe de Heisenberg $H(W)$ est le groupe d'ensemble sous-jacent
$W \oplus {\Bbb R}$, muni du produit
$$(w_1 , t_1) \cdot (w_2 , t_2) = (w_1 +w_2 , t_1 + t_2 + \frac12 B (w_1 , w_2) ).$$
Le groupe symplectique 
$$Sp(W) = \{ g \in GL(W) \; : \; B( g \cdot v , g \cdot w ) = B( v,w ) \}$$
agit comme groupe d'automorphismes de $H(W)$ via~:
$$g \cdot (w,t) = (g\cdot w , t).$$
Cette action est triviale sur le centre $Z= \{ (0,t) \} \cong {\Bbb R}$ de $H(W)$.

\paragraph{2.} Le th\'eor\`eme de Stone-von Neumann affirme qu'il existe une unique classe d'\'equivalence de repr\'esentation irr\'eductible unitaire $\rho$ de $H(W)$ de caract\`ere central
$t \mapsto e^{2i\pi t}$, {\it i.e.} telle que 
$$\rho ((0,t)) = e^{2i\pi t} \cdot {\rm Id}.$$
Nous notons $\rho^{\infty}$ la repr\'esentation lisse correspondante, {\it i.e.} sur les vecteurs $C^{\infty}$. Elle est \'egalement unique.
L'aspect int\'eressant de la th\'eorie des repr\'esentations du groupe d'Heisenberg n'est donc pas 
son dual unitaire mais plut\^ot l'ensemble des diff\'erentes r\'ealisations de la repr\'esentation $\rho$ (resp. $\rho^{\infty}$), les
{\it mod\`eles}. 

\paragraph{3. Mod\`eles de Schr\"odinger.} Une {\it polarisation compl\`ete} de $W$ est la donn\'ee d'une d\'ecomposition $W=X+Y$, o\`u $X$ et 
$Y$ sont deux sous-espaces totalement isotropes maximaux de $W$. La restriction de $B$ induit une dualit\'e 
$$\langle \cdot , \cdot \rangle : X \times Y \rightarrow {\Bbb R}.$$
(On a alors $B(x+y , x'+y') = \langle x , y' \rangle - \langle x' , y \rangle$, o\`u $x,x' \in X$ et $y,y' \in Y$.) 
Une telle polarisation donne lieu \`a une r\'ealisation 
de la repr\'esentation $\rho$ dans l'espace $L^2 (X)$ o\`u l'action de $H(W)$ est donn\'ee par 
$$\left( \rho ((x+y , t)) f \right) (x') = \exp \left( 2 i\pi (t - \langle x' - \frac12 x, y \rangle ) \right) f(x'-x),$$
o\`u $f \in L^2 (X)$, $x \in X$, $y\in Y$ et $t \in {\Bbb R}$. 

L'espace de la repr\'esentation $\rho^{\infty}$, autrement dit l'espace des vecteurs $C^{\infty}$ de $(\rho , L^2 (X))$, muni de la topologie $C^{\infty}$, est 
isomorphe \`a l'espace de Schwartz ${\cal S} (X)$ de $X$ muni de la topologie de Schwartz. L'alg\`ebre de Lie de $H(W)$ s'identifie \`a l'ensemble 
$W \oplus {\Bbb R}$ muni du crochet
$$[(w_1 , t_1) , (w_2 , t_2)] = (0, B(w_1 , w_2 )).$$
Elle est en particulier engendr\'ee par les vecteurs $e_1 , \ldots , e_N, f_1 , \ldots , f_N$ d'une base symplectique de $W$. Supposons alors 
que $(e_1 , \ldots , e_N)$ est une base de $X$ et $(f_1 , \ldots , f_N )$ la base duale de $Y$. Cette derni\`ere permet 
de d\'efinir les coordonn\'ees $(x_1 , \ldots , x_N)$ d'un \'el\'ement de $X$~: $x_j (x) = \langle x, f_j\rangle$ et
l'action infinit\'esimale de $\rho$ dans ${\cal S}(Y)$ est engendr\'ee par les op\'erateurs
\begin{eqnarray} \label{action1}
\begin{array}{l}
\rho (e_j ) = -\frac{\partial}{\partial x_j} , \\
\rho (f_j ) = -2i\pi x_j.
\end{array}
\end{eqnarray}

\paragraph{4. Mod\`eles de Fock.} Une structure complexe $J$ d\'efinie positive ({\it i.e.} la forme $Q(\cdot , \cdot ) = B (\cdot , J \cdot )$ est 
sym\'etrique d\'efinie positive) sur $W$ donne lieu \`a une r\'ealisation de la repr\'esentation $\rho$ dans l'espace de Hilbert 
${\cal F}$ compl\'et\'e de l'ensemble des fonctions $\phi$ $J$-holomorphes sur $W$ telles que
$$\int_W |\phi (w)|^2 e^{- 2\pi H(w,w)}  dw < +\infty,$$
o\`u $H$ est l'unique forme hermitienne sur $W$ de partie imaginaire $B$, \`a savoir $H=Q+iB$.
L'action de $H(W)$ dans ${\cal F}$ est donn\'ee par 
$$\left( \rho ((w_0,t)) \phi \right) (w) = \exp \left( 2i\pi t + \pi H(w_0, w -\frac12 w_0) \right) \phi (w -w_0),$$
o\`u $\phi \in {\cal F}$, $w,w_0 \in W$ et $t\in {\Bbb R}$.

Soit $(e_1 , \ldots , e_N, f_1 , \ldots , f_N)$ une base symplectique de $W$ telle que
$$J e_i = f_i \; \mbox{  et  } \; J f_i = -e_i \; \mbox{  pour  } 1 \leq i \leq N.$$
On peut d\'ecomposer $W\otimes {\Bbb C}$ en sous-espaces propres sous l'action de $J$~:
$$W\otimes {\Bbb C} = W' + W'',$$
o\`u $W'$ est associ\'e \`a la valeur propre $+i$ et $W''$ \`a la valeur propre $-i$. Le sous-espace $W'$ muni de sa structure complexe (multiplication
par $i$) s'identifie donc \`a l'espace $W$ muni de la structure $J$. Notons $(w_1 ' , \ldots , w_N ')$ et $(w_1 '' , \ldots , w_N '')$ les bases complexes respectives 
de $W'$ et $W''$ d\'efinies par
$$w_j' = e_j -i f_j  \; \mbox{  et  } \; w_j '' = e_j +if_j \; \mbox{ pour } 1\leq j \leq n .$$
La forme bilin\'eaire anti-sym\'etrique $B$ se prolonge naturellement \`a $W\otimes {\Bbb C}$ et induit, par restriction, une dualit\'e 
$$\langle \cdot , \cdot \rangle : W' \times W'' \rightarrow {\Bbb C}.$$
Notons $z_j$, pour $1\leq j \leq N$, les formes lin\'eaires sur $W'$ correspondantes aux \'el\'ements $w_j'' \in W''$ via cette dualit\'e~: 
$$z_j (w') = \langle w', w_j'' \rangle .$$

L'alg\`ebre des polyn\^omes holomorphes sur $W$ s'identifie naturellement \`a l'alg\`ebre sym\'etrique $S^* (W')^* (\cong S^* (W'')$ via la dualit\'e ci-dessus).
Cette alg\`ebre forme un sous-espace dense de vecteurs $C^{\infty}$ dans ${\cal F}$ stable sous l'action infinit\'esimale de $\rho^{\infty}$. 
Les $z_j$ d\'efinissent des coordonn\'ees complexes et permettent d'identifier $S^* (W')^*$ avec l'alg\`ebre ${\cal P} ({\Bbb C}^N)$ des polyn\^omes
en $z_1, \ldots , z_N$. Notons $\frac{\partial}{\partial z_j}$ la d\'erivation des polyn\^omes en $z_1,  \ldots , z_N$ d\'etermin\'ees par $\frac{\partial}{\partial z_j} (z_k ) =\delta_{jk}$.
L'action infinit\'esimale (complexifi\'ee) de $\rho$ dans ${\cal P} ({\Bbb C}^N)$ est alors engendr\'ee par les op\'erateurs
\begin{eqnarray} \label{action2}
\begin{array}{l}
\rho (w_j ' ) = -4\pi\frac{\partial}{\partial z_j} ,\\
\rho (w_j '' ) = z_j .
\end{array}
\end{eqnarray}

\paragraph{5.} Une structure complexe $J$ d\'efinie positive sur $W$ induit une polarisation compl\`ete $W=X+Y$ avec $Y=JX$. 
Il est alors possible de choisir une base symplectique $(e_1 , \ldots , e_N, f_1 , \ldots , f_N)$ de $W$ v\'erifiant les conditions des paragraphes (3 et 4) pr\'ec\'edents. 
On peut donc chercher \`a d\'ecrire l'op\'erateur d'entrelacement entre les mod\`eles de Fock et de Schr\"odinger correspondants. Remarquons qu'en
conservant les notations pr\'ec\'edentes, $W' = X \otimes {\Bbb C}$, $W'' = Y\otimes {\Bbb C}$ et $x_j = {\rm Re}(z_j)$ pour $j=1 , \ldots , N$.

L'action de $W'$ sur ${\cal S}(X)$, induite par (\ref{action1}) est d\'etermin\'ee par 
$$\rho (w_j ') = - \frac{\partial}{\partial x_j} - 2\pi x_j \; \mbox{  pour } j=1, \ldots , N.$$
La gaussienne 
\begin{eqnarray} \label{gauss}
\varphi_0 (x ) = \exp ( - \pi H(x,x)) \in {\cal S}(X)
\end{eqnarray} 
est donc annul\'ee par l'action de $W'$. Mais il est imm\'ediat dans le mod\`ele de Fock que le sous-espace annul\'e par l'action de $W'$ est de dimension un
et constitu\'e des fonctions constantes. (Une fonction holomorphe annul\'ee par tous les op\'erateurs $\frac{\partial}{\partial z_j}$ est constante.)  
Il existe donc un unique op\'erateur d'entrelacement entre le mod\`ele de Fock et le mod\`ele de Schr\"odinger qui envoie la fonction constante \'egale \`a 
$1$ sur la gaussienne $\varphi_0$. 
Notons ${\rm S}(X)$ l'image de ${\cal P}({\Bbb C}^N)$ dans ${\cal S}(X)$ via cet op\'erateur d'entrelacement. Le sous-espace
${\rm S}(X) \subset {\cal S}(X)$ est stable sous l'action infinit\'esimale de $\rho^{\infty}$; il co\"incide avec le sous-espace $\varphi_0 {\Bbb R} [x_1 , \ldots , x_N]$.
Les op\'erateurs correspondants respectivement \`a la d\'erivation $\partial/\partial z_j$ et \`a la multiplication par $z_j$ dans ${\cal P}[{\Bbb C}^N]$ sont respectivement (et \`a des constantes multiplicatives pr\`es)~:
$$A_j^+ = \frac{\partial}{\partial x_j} + 2\pi x_j  \; \mbox{  et  } \; A_j^- = \frac{\partial}{\partial x_j} - 2\pi x_j .$$
L'espace ${\rm S}(X)$ est donc invariant par ces deux op\'erateurs qui v\'erifient~:
\begin{enumerate} 
\item $A_j^+ \varphi_0 = 0$,
\item $[A_j^+ , A_i^+] = [A_j^- , A_i^- ]=0$,
\item $[A_i^+ , A_j^-] = -4\pi \delta_{ij} {\rm Id}$ et,
\item si $H_j := A_j^+ A_j^- + A_j^- A_j^+ = 2 \left( \frac{\partial^2}{\partial y_j^2} -4\pi^2 y_j^2 \right)$, 
$$[H_j , A_i^+ ] = 8\pi \delta_{ij} A_i^+ \; \mbox{  et  } \; [H_j , A_i^- ] = -8\pi \delta_{ij} A_i^- .$$
\end{enumerate}
On v\'erifie par ailleurs que chaque op\'erateur $H_j$ est auto-adjoint dans ${\rm S}(X)$ et que $H_j \varphi_0 = -4\pi \varphi_0$, pour $j=1, \ldots , N$.

Si $m=(m_1 , \ldots , m_N) \in {\Bbb N}^N$, posons
$$\varphi_m = (m_1 ! \ldots m_N !)^{-1/2} (A_1^- )^{m_1} \ldots (A_N^- )^{m_N } \varphi_0 .$$
Il d\'ecoule facilement des formules ci-dessus que la famille $\{\varphi_m \}_{m \in {\Bbb N}^N}$ forme une base orthogonale de $L^2 (X)$. 

\paragraph{6. Mod\`eles mixtes I.} Une d\'ecomposition {\it mixte} de $W$ est une d\'ecomposition
$$W= X+W_0 +Y,$$ 
o\`u $X$ et $Y$ sont deux sous-espaces isotropes en dualit\'e
pour la restriction de $B$ et o\`u $W_0$, le suppl\'ementaire orthogonale de $X+Y$ dans $W$, est muni d'une structure complexe positive $J_0$ (pour la restriction de 
$B$ \`a $W_0$). Une d\'ecomposition mixte donne lieu \`a une r\'ealisation de la repr\'esentation $\rho$ dans l'espace $L^2 (X, {\cal F}_0 )$
des fonctions sur $X$ \`a valeurs dans ${\cal F}_0$ et de norme de carr\'e int\'egrable., o\`u ${\cal F}_0$ d\'esigne
le mod\`ele de Fock associ\'e \`a $(W_0 , J_0)$. Le sous-groupe $H(W_0)$ est distingu\'e dans le groupe $H(W)$;  il agit sur $L^2 (X , {\cal F}_0 )$ via sa repr\'esentation
$\rho_0$ dans ${\cal F}_0$~:
$$\left( \rho ((w_0 , t)) f \right) (x) = \rho ((w_0 , t)) (f(x)).$$
L'action de $H(W)$ sur $L^2 (X , {\cal F}_0 )$ est finalement compl\`etement d\'ecrite par~:
$$\left( \rho ((x+y , t)) f \right) (x') = \exp \left( 2i\pi (t - \langle x' - \frac12 x, y \rangle ) \right) f(x'-x),$$
o\`u $f \in L^2 (X, {\cal F}_0 )$, $x \in X$, $y\in Y$ et $t \in {\Bbb R}$. 

Il est imm\'ediat que l'ensemble des vecteurs $C^{\infty}$ de $\rho$ dans cette r\'ealisation est l'espace de Schwartz ${\cal S} (X , {\cal F}^{\infty})$, de mani\`ere analogue
\`a la d\'efinition de ${\rm S}(X)$, il correspond aux polyn\^omes de la r\'ealisation de Fock, le sous-espace dense ${\rm S}(X , {\cal P} ({\Bbb C}^{N_0}))$ de 
${\cal S} (X , {\cal F}^{\infty})$.

\paragraph{7. La repr\'esentation de l'oscillateur.} 
Il existe un unique rev\^etement non trivial \`a deux feuillets $\widetilde{Sp}(W)$ du groupe $Sp(W)$, appel\'e {\it groupe m\'etaplectique}. Notons
$\widetilde{g} \mapsto g$ la projection de rev\^etement. Le th\'eor\`eme de Shale-Weil affirme qu'il existe une unique classe d'\'equivalence de repr\'esentation unitaire
$\omega$ du groupe $\widetilde{Sp}(W)$ telle que 
$$\omega (\widetilde{g} ) \rho ((w,t)) \omega (\widetilde{g}^{-1}) = \rho (g\cdot (w,t)).$$
Il correspond \`a $\omega$ la repr\'esentation lisse $\omega^{\infty}$. Les vecteurs lisses de $\rho$ et $\omega$ co\"{\i}ncident (dans n'importe quel r\'ealisation commune); on peut
donc r\'ealiser les repr\'esentations $\rho^{\infty}$ et $\omega^{\infty}$ dans un m\^eme espace. 

L'alg\`ebre de Lie du groupe m\'etaplectique co\"{\i}ncide avec l'alg\`ebre de Lie $\mathfrak{sp}(W)$ du groupe symplectique $Sp(W)$. Elle contient un tore compact
$i\mathfrak{t}_0 \subset \mathfrak{sp}(W)$ dont l'action infinit\'esimale sur ${\cal S}(X)$ est engendr\'ee par les op\'erateurs 
$iH_j$, $j=1, \ldots ,N$. Il en d\'ecoule facilement, que l'espace des vecteurs $C^{\infty}$ de la restriction de $\omega$ au tore $T= \exp (i\mathfrak{t}_0 ) \subset \widetilde{Sp}(W)$
co\"{\i}ncide avec l'espace des vecteurs $C^{\infty}$ de $\omega$.

Le sous-groupe $U(N) \subset Sp(W)$ constitu\'e des \'el\'ements qui centralise $J$ est un sous-groupe compact maximal; il s'identifie au groupe unitaire de la
forme hermitienne $H$. L'espace des vecteurs $\widetilde{U(N)}$-finis de la repr\'esentation $\omega$ dans un mod\`ele de Fock (resp. de Schr\"odinger) est 
${\cal P} ({\Bbb C}^N )$ (resp. ${\rm S}(X)$). L'action de l'alg\`ebre de Lie complexifi\'ee $\mathfrak{sp} (W)$ de $Sp(W)$ sur ${\cal P}({\Bbb C}^N)$ peut \^etre 
d\'ecrite par les op\'erateurs suivants~:
\begin{eqnarray} \label{omegainf}
\omega (\mathfrak{sp} (W) ) = \mathfrak{sp} (W)^{(1,1)} \oplus \mathfrak{sp} (W)^{(2,0)} \oplus \mathfrak{sp} (W)^{(0,2)} ,
\end{eqnarray}
o\`u
\begin{eqnarray} \label{omegainf2}
\begin{array}{ccl}
\mathfrak{sp} (W)^{(1,1)} & = & {\rm Vect} \left\{ -i \left( z_k \frac{\partial}{\partial z_j} + \frac{\partial }{\partial z_j} z_k \right) \right\}, \\
\mathfrak{sp} (W)^{(2,0)} & = & {\rm Vect} \{ i z_j z_k \}, \\
\mathfrak{sp} (W)^{(0,2)} & = & {\rm Vect} \left\{ 4i \frac{\partial^2 }{\partial z_j \partial z_k} \right\} .
\end{array}
\end{eqnarray} 
La d\'ecomposition (\ref{omegainf}) correspond \`a la d\'ecomposition de Cartan
\begin{eqnarray} \label{cartan}
\mathfrak{sp} (W) = \mathfrak{k} \oplus \mathfrak{p}_+ \oplus \mathfrak{p}_-, 
\end{eqnarray}
o\`u $\mathfrak{sp} (W)^{(1,1)} \cong \omega (\mathfrak{k})$, $\mathfrak{sp} (W)^{(2,0)} \cong \omega (\mathfrak{p}^+ )$ et $\mathfrak{sp} (W)^{(0,2)} \cong \omega (\mathfrak{p}^- )$.
Si ${\cal P} ({\Bbb C}^N) = \sum_{d \geq 0} {\cal P}^d ({\Bbb C}^N)$ est la graduation par le degr\'e de l'alg\`ebre des polyn\^omes ${\cal P}({\Bbb C}^N )$, il 
est imm\'ediat que $\mathfrak{sp}(W)^{(i,j)}$ envoie ${\cal P}^d ({\Bbb C}^N)$ sur ${\cal P}^{d+i-j} ({\Bbb C}^N)$.

\paragraph{8. Paires r\'eductives duales.} \`A la suite de Howe \cite{Howe}, on dit qu'un couple $(G,G')$ de sous-groupes de $Sp(W)$ forme une {\it paire r\'eductive duale} si 
\begin{enumerate}
\item $G$ et $G'$ agissent de mani\`eres absolument r\'eductibles sur $W$; et
\item $G$ est le centralisateur de $G'$ dans $Sp(W)$ et vice versa. 
\end{enumerate}
Nous ne consid\'erons ici que des paires r\'eductives duales {\it de type I}, {\it i.e.} telles que l'action de $G\cdot G'$ soit irr\'eductible sur $W$; celles sont associ\'ees \`a une alg\`ebre \`a division
$D$ sur ${\Bbb R}$ muni d'une involution. Dans la suite nous supposerons toujours $D = {\Bbb R}$ ou ${\Bbb C}$ et l'involution respectivement triviale ou la conjugaison complexe usuelle. Il existe alors
\begin{enumerate}
\item deux $D$-modules $V$ et $V'$ munis de formes sesquilin\'eaires respectivement not\'ees $(\cdot, \cdot)$ et $(\cdot , \cdot)'$; l'une hermitienne et l'autre 
anti-hermitienne; de telle mani\`ere que
\item $G$ et $G'$ soient les groupes d'isom\'etries respectifs de $(\cdot, \cdot)$ et $(\cdot , \cdot)'$, et
\item $W=V\otimes_D V'$, $B = {\rm tr}_{D/{\Bbb R}} \left( (,)\otimes \overline{(,)^{' }} \right)$.
\end{enumerate}

Dans la suite nous notons 
$$m = \dim_{D} V , \; m' =\dim_D V' , d = \dim_{\Bbb R} D,$$
et 
$$\varepsilon ' = \left\{
\begin{array}{cl}
1 & \mbox{ si } (\cdot , \cdot )' \mbox{ est } \mbox{hermitienne}, \\
-1 & \mbox{ sinon.}
\end{array} \right.$$

Si $H$ est un sous-groupe de $Sp(W)$, nous notons $\widetilde{H}$ sa pr\'eimage dans $\widetilde{Sp}(W)$. 
Commen\c{c}ons par traiter l'exemple de la paire duale $(U(r),U(p,q))$.

\paragraph{9. La paire $(U(r) , U(p,q))$.} (cf. Kashiwara-Vergne \cite{KashiwaraVergne}.)
Pour $p+q=n$, $p\geq q \geq 0$, posons
$$I_{p,q} = \left( 
\begin{array}{cc}
1_p & \\
 & -1_q 
\end{array} \right),$$
o\`u $1_p$ d\'esigne la matrice identit\'e de taille $p\times p$. Si $g \in M_{n} ({\Bbb C})$, posons $g^*$ \'egal au conjugu\'e de la transpos\'ee de la matrice $g$. 
Soit $(\cdot , \cdot )$ la forme hermitienne sur $V={\Bbb C}^n$ de matrice $I_{p,q}$. Le groupe $G=U(p,q)$ des isom\'etries de $(\cdot , \cdot)$ est 
le groupe de tous les $g \in M_n ({\Bbb C})$ tels que 
$$g^* I_{p,q} g = I_{p,q}.$$

Cons\'erons l'espace hermitien $(V' , (\cdot , \cdot )')$, o\`u $V' = {\Bbb C}^r$ et $(\cdot , \cdot )'$ de matrice $1_r$ et notons $G' = U(r)$ des isom\'etries
de $(\cdot , \cdot )'$. L'espace $W= V \otimes_{{\Bbb C}} V'$ est naturellement muni d'une forme symplectique $B = - 2 {\rm Im} \left( (,) \otimes \overline{(,)'} \right)$ et
le couple $(U(p,q) , U(r))$ forme un paire r\'eductive duale dans $Sp(W)$.

Nous allons consid\'erer un mod\`ele de Fock de la repr\'esentation $\omega$ de l'oscillateur du groupe $Sp(W)=Sp(2rn)$ restreinte \`a la paire $(U(p,q) , U(r))$. Le mod\`ele de 
Fock que nous consid\'erons est associ\'e \`a la structure complexe d\'efinie positive
$$J_{0} = J \otimes I_{p,q},$$
o\`u $J$ est la structure complexe usuelle ``multiplication par $i$'' sur ${\Bbb C}^r$ et $I_{p,q}$ est l\`a pour tordre $J$ en une structure complexe d\'efinie positive.

Dans le mod\`ele de Fock associ\'e, l'espace des vecteurs $\widetilde{U(rn)}$-finis s'identifie avec l'anneau des polyn\^omes ${\cal P} (M_{p\times r} \oplus M_{q\times r} )$. 
C'est aussi l'espace des vecteurs $(U(r) \cdot (U(p) \times U(q)))$-finis de la restriction de $\omega$ au groupe $G\cdot G'$ (le groupe $U(r) \cdot (U(p) \times U(q))$ contient un tore compact de
$Sp(2rn)$). Notons enfin que bien que $\omega$ ne soit pas une vraie repr\'esentation du groupe $Sp(W)$ mais seulement de son rev\^etement m\'etaplectique, elle se 
restreint en une vraie repr\'esentation du groupe $U(r) \cdot U(p,q)$.

Les sous-groupes compacts maximaux des groupes 
$U(r)$ et $U(p,q)$ se complexifient respectivement en $GL(r)$ et $GL(p) \times GL(q)$; ils agissent sur ${\cal P} (M_{p\times r} \oplus M_{q\times r} )$, via $\omega$, de la mani\`ere suivante~:
\begin{eqnarray} \label{actP}
(g,h_1 , h_2 ) \cdot P(X,Y) = (\det h_2 )^r P(h_1^{-1} X g , {}^t h_2 Y {}^t g^{-1} ),
\end{eqnarray}
o\`u  $X \in M_{p \times r}$, $Y\in M_{q \times r}$, $ g \in GL(r)$, $h_1 \in GL(p)$ et $h_2 \in GL(q)$. 

Notons $E={\Bbb C}^p$ (resp. $F={\Bbb C}^q$, $G={\Bbb C}^r$). La repr\'esentation de $GL(r) \times GL(p) \times GL(q)$ ci-dessus correspond alors \`a la repr\'esentation
de $GL(G) \times GL(E) \times GL(F)$ dans 
$$D_r (F) \otimes {\rm Sym} [E^* \otimes G] \otimes {\rm Sym} [F \otimes G^*],$$
o\`u nous notons $D_l = (\bigwedge^k E)^{\otimes l}$ si $l\geq 0$ et $D_{-l}$ est la repr\'esentation duale $D_l^*$.

Choisissons comme sous-alg\`ebre de Borel dans $\mathfrak{u}(k) \times \mathfrak{u} (r)$ l'alg\`ebre
des matrices qui sont triangulaires sup\'erieures sur $E$ et triangulaires inf\'erieures sur $F$ par rapport aux bases canoniques de $E$ et $F$.
Il est alors bien connu (cf. \cite{Fulton}) qu'\`a chaque partition $\lambda$ de longueur $l(\lambda ) \leq p$, il correspond une repr\'esentation irr\'eductible finie $E^{\lambda}$ du groupe $GL(E)$. 
Toutes les repr\'esentations du groupe $GL(E)$ sont d'ailleurs obtenues en formant le produit tensoriel 
$E^{\lambda} \otimes D_l$ ($l\in {\Bbb Z}$) et sont deux \`a deux non isomorphes. 

L'int\'er\^et du formalisme des paires r\'eductives duales provient de la d\'ecomposition en irr\'eductible de la restriction \`a $G \cdot G'$ de la repr\'esentation de l'oscillateur $\omega$.
Le cas $q=0$ est d\'ej\`a important. On sait en effet d\'ecomposer en irr\'eductibles la repr\'esentation de 
$GL(G) \times GL(E)$ dans ${\rm Sym} [E^* \otimes G]$ (cf. \cite{Fulton})~:
$${\rm Sym} [E^* \otimes G ] = \bigoplus_{\lambda } (E^{\lambda})^* \otimes G^{\lambda}.$$ 
(Ici $E^{\lambda}$ (resp. $G^{\lambda}$) est trivial si $l(\lambda) > p$ (resp. $r$).)
Finalement, la repr\'esentation $\omega$ qui, restreinte au groupe $\widetilde{U} (r) \cdot \widetilde{U} (p)$, d\'efinit une vraie 
repr\'esentation du groupe $U(r) \cdot U(p)$, se d\'ecompose en irr\'eductibles et cette d\'ecomposition d\'efinit une correspondance bijective entre certaines repr\'esentations des groupes 
$U(r)$ et $U(p)$, la {\it correspondance th\^eta}, que nous d\'ecrivons dans la proposition suivante.

\begin{prop}
Soit $\lambda$ un diagramme de Young. La repr\'esentation $G^{\lambda}$ de $U(r)$ intervient dans la correspondance th\^eta si et seulement si 
sa longueur $l(\lambda ) \leq \min (p,r)$; auquel cas la repr\'esentation correspondante de $U(p)$ est la repr\'esentation $(E^{\lambda})^*$.
\end{prop}

Il reste plus g\'en\'eralement vrai que la restriction de la repr\'esentation $\omega$ de $\widetilde{Sp} (2r(p+q))$ \`a $\widetilde{U}(r) \cdot \widetilde{U}(p,q)$ d\'efinit une vraie 
repr\'esentation du groupe $U(r) \cdot U(p,q)$. Notons 
$$\mathfrak{u}(p,q)^{(i,j)} = \mathfrak{sp}(2r(p+q))^{(i,j)} \cap \omega (\mathfrak{u}(p,q)).$$
On obtient alors la d\'ecomposition 
\begin{eqnarray} \label{dec}
\omega (\mathfrak{u}(p,q)) = \mathfrak{u}(p,q)^{(1,1)} \oplus \mathfrak{u}(p,q)^{(2,0)} \oplus \mathfrak{u}(p,q)^{(0,2)}
\end{eqnarray}
qui n'est autre que la d\'ecomposition de Cartan complexifi\'ee
\begin{eqnarray} \label{cartan}
\mathfrak{k} \oplus \mathfrak{p}_+ \oplus \mathfrak{p}_- , 
\end{eqnarray}
o\`u $\mathfrak{u}(p,q)^{(1,1)} = \omega (\mathfrak{k})$, $\mathfrak{u}(p,q)^{(2,0)} = \omega (\mathfrak{p}^+)$ et $\mathfrak{u}(p,q)^{(0,2)} = \omega (\mathfrak{p}^-)$. En particulier,
$\mathfrak{k}$ a un centre de dimension un et $\mathfrak{p}^{\pm}$ sont les $\pm i$-espaces propres de ce centre. 

La loi de branchement de $\omega_{| U(r) \cdot U(p,q)}$ est alors de la forme
\begin{eqnarray} \label{B1}
{\cal P} (M_{p\times r} \oplus M_{q\times r} )_{|U(r) \cdot U(p,q)} = \bigoplus_{\tau \in S} \tau \otimes V_{\tau '} ,
\end{eqnarray}
o\`u $S$ est un sous-ensemble du dual unitaire de $U(r)$. Les repr\'esentations $V_{\tau '}$ de $U(p,q)$ sont irr\'eductibles et {\it holomorphes}, {\it i.e.} il existe un vecteur non nul $v \in V_{\tau '}$ tel que 
$\tau ' (\mathfrak{p}^- ) \cdot v =0$. L'aspect important de cette d\'ecomposition est l'unicit\'e de cette correspondance, {\it i.e.} une repr\'esentation $\tau$ de $U(r)$ appara\^{\i}t une seule fois
et d\'etermine une unique repr\'esentation de $U(p,q)$, c'est la {\it correspondance th\^eta}. 

Soit ${\cal H} = \ker \mathfrak{u}(p,q)^{(0,2)} \subset {\cal P} (M_{p\times r} \oplus M_{q\times r} )$ l'espace des polyn\^omes harmoniques, {\it i.e.} l'espace des polyn\^omes $P(X,Y)$ tels que
\begin{eqnarray} \label{eqdiff}
(\Delta_{ij} P)(X,Y) =0 \; \mbox{ pour } 1 \leq i \leq p, \; 1 \leq j \leq q,
\end{eqnarray}
o\`u 
$$\Delta_{ij} = \sum_{\nu =1}^r \frac{\partial^2}{\partial X_{i\nu} \partial Y_{j\nu }} , \; X \in M_{p\times r} , \; Y \in M_{q \times r} .$$

\begin{thm}[Kashiwara-Vergne] \label{KV}
L'espace ${\cal H}$ est un $U(r) \times U(p) \times U(q)$-module qui admet une d\'ecomposition en irr\'eductibles sans multiplicit\'e~:
\begin{eqnarray} \label{decH}
{\cal H} = \bigoplus_{\tau \in S} \tau \otimes \ker \tau ' (\mathfrak{u}(p,q)^{(0,2)}).
\end{eqnarray}
Et la correspondance th\^eta est d\'ecrite par~:
\begin{eqnarray} \label{thetacorres}
\begin{array}{ccl}
{\cal P} (M_{p\times r} \oplus M_{q\times r} ) & = & {\cal H} \cdot {\cal S} (\mathfrak{u}(p,q)^{(2,0)}) ,\\
& = & \bigoplus_{\tau \in S} \tau \otimes \left\{ {\cal S} (\mathfrak{u}(p,q)^{(2,0)}) \cdot \ker \tau ' (\mathfrak{u}(p,q)^{(0,2)})\right\} ,\\
& = & \bigoplus_{\lambda , \mu } \tau (\lambda , \mu) \otimes  M(\lambda , \mu ),
\end{array}
\end{eqnarray}
o\`u la derni\`ere somme porte sur l'ensemble des partitions $\lambda$ et $\mu$ telles que $l(\lambda ) \leq p$, $l(\mu ) \leq q$ et $l(\lambda ) + l(\mu ) \leq r$, la repr\'esentation 
$\tau(\lambda , \mu)$ de $U(r)$ est l'unique repr\'esentation de plus haut poids $(\lambda , 0 , -\mu )$ et $M(\lambda , \mu)$ est un $U(p,q)$-module holomorphe de plus haut poids d\'etermin\'e
par son plus bas $K$-type
$$(E^{\lambda})^* \otimes F^{\mu} \otimes D_r (F).$$
\end{thm}

\'Etant donn\'es $X\in M_{p \times r}$ et $Y \in M_{q \times r}$, posons~: 
$$\Delta_j (X) = \det \left(
\begin{array}{ccc} 
X_{11} & \ldots & X_{1 j} \\
\vdots & & \vdots \\
X_{j1} & \ldots & X_{jj} 
\end{array} \right) \;  \;  (0\leq j \leq r,p) , $$
$$\widetilde{\Delta}_j (Y) = \det \left(
\begin{array}{ccc} 
Y_{q-j+1,r-j+1} & \ldots & Y_{q-j+1 , r} \\
\vdots & & \vdots \\
Y_{q,r-j+1} & \ldots & Y_{qr} 
\end{array} \right) \;  \;  (0\leq j \leq r,q) .$$

\begin{prop} \label{KV2} 
Le polyn\^ome 
$$P_{\lambda, \mu} (X,Y) = \Delta_1 (X)^{\lambda_1 - \lambda_2} \ldots \Delta_{l(\lambda)} (X)^{\lambda_{l(\lambda)}} \cdot \widetilde{\Delta}_1 (Y)^{\mu_1 - \mu_2} \ldots \widetilde{\Delta}_{l(\lambda)} (Y)^{\mu_{l(\mu)}}$$
dans ${\cal P}(M_{p\times r} \oplus M_{q\times r})$ correspond, dans le mod\`ele de la repr\'esentation de Weil d\'ecrit ci-dessus, \`a un vecteur de plus haut poids pour le groupe 
$U(p) \times U(q) \times U(r)$. Sous l'action de $U(r)$, il engendre la repr\'esentation $\tau (\lambda , \mu )$, sous l'action de $U(p,q)$, il engendre la repr\'esentation $M(\lambda , \mu)$.
\end{prop}

\paragraph{13. Paires en balance.} Deux paires r\'eductives duales $(G,H')$ et $(H,G')$ dans un m\^eme groupe symplectique $Sp(W)$ sont 
dites {\it en balance} si 
$$H \subset G \mbox{ et } H' \subset G'.$$

Le diagramme qui motive cette terminologie, due \`a Kudla \cite{Kudla}, est le suivant~:
$$
\begin{array}{ccc}
G & & G' \\
|  & \times & | \\
H &  & H' 
\end{array}
$$
o\`u les droites diagonales relient les membres d'une m\^eme paire duale et les droites verticales correspondent aux inclusions de groupes.

Les repr\'esentations (de Weil) de $G\cdot H'$ et $H \cdot G'$ peuvent \^etre r\'ealis\'ees dans un m\^eme espace et se restreignent en une m\^eme
repr\'esentation de $H \cdot H'$ sans qu'aucune repr\'esentation de $G \cdot G'$ n'existe. C'est une mani\`ere commode de comprendre d'o\`u proviennent
certaines identit\'es ``miraculeuses'' entre s\'eries th\^eta. Nous utilisons des paires en balances dans la sections suivantes pour construire les
classes du th\'eor\`eme \ref{T1} par cup-produit et restriction \`a partir de classes explicites d\'eduites du th\'eor\`eme de Kashiwara et Vergne.

\section{Construction de classes de Fock/Schwartz}

\subsection{Classes (anti-)holomorphes dans les groupes unitaires}

Comme au \S 9 de la section pr\'ec\'edente, dans ce paragraphe nous notons $(\cdot , \cdot )$ la forme hermitienne sur $V={\Bbb C}^n$ de matrice $I_{p,q}$, $G=U(p,q)$ le groupe des isom\'etries de $(\cdot , \cdot)$,
$(V' , (\cdot , \cdot )')$ l'espace hermitien $V' = {\Bbb C}^r$ de forme hermitienne $(\cdot , \cdot )'$ associ\'ee \`a la matrice $1_r$ et $G' = U(r)$ le groupe des isom\'etries
de $(\cdot , \cdot )'$. Rappelons que nous avons d\'ecrit la repr\'esentation de l'oscillateur $\omega$ restreinte au groupe $G \cdot G'$ dans un mod\`ele de 
Fock que nous notons ici ${\cal F}$.

\medskip

Soit $D(=D^+)$ l'espace sym\'etrique associ\'e \`a $G$~:
$$D = G/K =U(p,q) /U(p) \times U(q).$$
Soit $\mathfrak{g}_0 = {\rm Lie} (G)$ l'alg\`ebre de Lie de $G$ et soit $\mathfrak{k}_0$ celle 
de $K=U(p) \times U(q)$. On a la d\'ecomposition
\begin{eqnarray} 
\mathfrak{g} = \mathfrak{k} + \mathfrak{p}_+ + \mathfrak{p}_-,
\end{eqnarray}
et l'espace des formes diff\'erentielles sur $D$ de type $(a,b)$ est
\begin{eqnarray}
\Omega^{a,b} (D ) \cong \left[ C^{\infty} (G) \otimes (\bigwedge {}^a (\mathfrak{p}_+ )^* \otimes 
\bigwedge {}^b (\mathfrak{p}_- )^* ) \right]^K .
\end{eqnarray}

Dans cette section nous construisons un \'el\'ement
\begin{eqnarray} \label{phi}
\psi \in [{\cal F}^{\infty} \otimes \Omega^{0,rq} (D ) ]^G \cong [{\cal F}^{\infty} \otimes \bigwedge {}^{0, rq} \mathfrak{p}^* ]^K.
\end{eqnarray}
Nous montrons que la forme diff\'erentielle ainsi d\'efinie est ferm\'ee, {\it i.e.}
\begin{eqnarray} \label{d=0}
d \psi = 0,
\end{eqnarray}
o\`u $d$ est la diff\'erentielle du complexe d\'efini par (\ref{phi}) qui co\"{\i}ncide avec le complexe
$${\rm Hom}_K (\bigwedge {}^{a , b} \mathfrak{p} , {\cal F}^{\infty} )$$
calculant la $(\mathfrak{g} , K)$-cohomologie du $(\mathfrak{g} , K)$-module ${\cal F}^{\infty}$. Nous produisons 
ainsi une classe de cohomologie $[\psi ] \in H^{rq} (\mathfrak{g} , K ; {\cal F}^{\infty} )$ dont nous v\'erifierons qu'elle est non nulle.

Notre construction est en fait plus g\'en\'erale, elle concerne d'autres types de classes de cohomologie. 
\'Etant donn\'e un couple d'entiers $(a,b)$ avec $0 \leq a \leq p$, $0 \leq b \leq q$ et $a+b=r$, notons 
$$\nu = ( \underbrace{q, \ldots , q}_{a \ {\rm fois}} , \underbrace{b , \ldots , b}_{p-a \ {\rm fois}} )$$
la partition de diagramme de Young~:
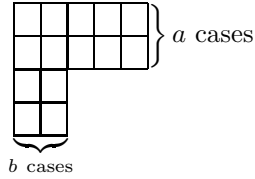
$$
\begin{array}{l} 
\left. \hspace{0,035cm}
\begin{array}{|c|c|c|c|c|} \hline
 & & & &\\ \hline
 & & & &\\ \hline
\end{array} 
\right\}  a \ {\rm cases} \\
\underbrace{
\begin{array}{|c|c|} \hline
 &  \\ \hline
 &  \\ \hline
\end{array}}_{b \ {\rm cases}}
\end{array} 
$$ 
(Ici $p=4$ et $q=5$.)

Il existe (cf. \cite{Lefschetz}) un unique $U(p,q)$-module anti-holomorphe, $\overline{A(\nu)}$, de plus haut poids dont le plus bas $K$-type est  
$$\overline{V(\nu )} = (E^{\nu })^* \otimes F^{\nu^*}.$$
(Ici on note encore $E={\Bbb C}^p$ et $F={\Bbb C}^q$.)

\noindent
Ce module est un module cohomologique anti-holomorphe; r\'eciproquement on obtient de cette mani\`ere et en faisant varier $r$ tous les modules
cohomologiques anti-holomorphes.

Le module $\overline{A(\nu)}$ n'intervient pas directement dans la correspondance th\^eta mais il d\'ecoule du th\'eor\`eme \ref{KV} que le module 
en correspondance avec 
la repr\'esentation de $U(r)$ de plus haut poids 
$$(\underbrace{q-b , \ldots , q-b}_{a \ {\rm fois}} , 0 , \ldots , 0, \underbrace{a-p , \ldots , a-p}_{b \ {\rm fois}} ) $$
co\"{\i}ncide, en restriction au groupe $SU(p,q)$, avec $\overline{A(\nu)}$.
(Ces deux modules ne diff\`erent, en effet que d'un caract\`ere central $(\det )^b$.)

\noindent
Par abus de notation, nous noterons $\overline{A(\nu)}$ le module $M(\lambda , \mu)$ o\`u 
$$\lambda =((q-b)^a ) \  \mbox{ et } \  \mu = ((p-a )^b).$$
Ce module est cohomologique pour le groupe $SU(p,q)$.

Le plus bas $K$-type du module $\overline{A(\nu)}$ est r\'ealis\'e dans l'espace des polyn\^omes ${\cal P} (M_{p\times r} \oplus M_{q\times r})$ pour
l'action (\ref{actP}). Il correspond au plongement 
$$(E^{\lambda})^* \otimes F^{\mu} \otimes D_r (F) \hookrightarrow D_{r} (F) \otimes \left( (E^{\lambda})^* \otimes G^{\lambda} \right) \otimes \left( F^{\mu} \otimes (G^{\mu})^*\right) ,$$
dans le sous-$U(r) \times U(p) \times U(q)$-module irr\'eductible de 
$$D_r (F) \otimes {\rm Sym} [E^* \otimes G] \otimes {\rm Sym} [F \otimes G^*] \cong {\cal P} (M_{p\times r} \oplus M_{q\times r})$$ 
engendr\'e par le vecteur de plus haut poids 
$$P_{\lambda , \mu } (X,Y) = \Delta_a (X)^{q-b} \cdot \widetilde{\Delta}_b (Y)^{p-a}   \in {\cal P} (M_{p\times r} \oplus M_{q\times r}).$$
Chacun des plongements $(E^{\lambda})^* \hookrightarrow {\Bbb C} [M_{p\times r}]$, $F^{\mu} \hookrightarrow {\Bbb C} [M_{q \times r}]$ correspond aux r\'ealisations classiques, {\it cf.} Fulton \cite[pp. 109--111] {Fulton}, 
de $(E^{\lambda})^*$ (resp. $F^{\mu}$) dans des espaces de polyn\^omes. On peut expliciter ces plongements sur des bases.

L'espace $E^{\lambda}$ s'obtient comme quotient de $(\bigwedge^a E) \otimes \ldots \otimes (\bigwedge^a E)$ ($q-b$ fois) par le sous-espace $Q^{\lambda} (E)$ engendr\'e par les 
\'el\'ements de la forme $\wedge \mathbf{v} - \sum \wedge \mathbf{w}$, o\`u l'on somme sur tous les \'el\'ements $\mathbf{w}$ obtenus \`a partir de $\mathbf{v}$ par un \'echange entre deux colonnes donn\'ees, d'un 
sous-ensemble donn\'e de boites dans la colonne de droite choisie ({\it cf.} \cite[\S 8.1]{Fulton} pour plus de d\'etails).
Soient $(e_1 , \ldots , e_p)$ et $(f_1 , \ldots , f_q)$ les bases canoniques respectives de $E$ et $F$. Les projet\'es des \'el\'ements 
$$e_T := \left( \bigwedge_{k=1}^a e_{i_{k,1}} \right) \otimes \ldots \otimes  \left( \bigwedge_{k=1}^a e_{i_{k,q-b}} \right),$$
o\`u $T$ est le {\it tableau} de Young de diagramme $\lambda$ et dont la case de coordonn\'ees $(k,l)$ est num\'erot\'ee $i_{k,l} \in \{ 1 , \ldots , p \}$, forment
une base de $E^{\lambda}$. Remarquons que l'hypoth\`ese que $T$ est un tableau signifie ici simplement que
$$i_{1,l} < \ldots < i_{a,l} \ \mbox{ pour tout } l= 1 , \ldots , q-b$$
et
$$i_{k,1} \leq \ldots \leq i_{k, q-b} \ \mbox{ pour tout } k =  1 , \ldots , a .$$
De la m\^eme mani\`ere, les projet\'es des \'el\'ements 
$$f_U := \left( \bigwedge_{k=1}^b f_{i_{k,1}} \right) \otimes \ldots \otimes  \left( \bigwedge_{k=1}^b f_{i_{k,p-a}} \right),$$
o\`u $U$ est le {\it tableau} de Young de diagramme $\mu$ et dont la case de coordonn\'ees $(k,l)$ est num\'erot\'ee $i_{k,l} \in \{ 1 , \ldots , q \}$, forment
une base de $F^{\mu}$.

Le plongement $(E^{\lambda})^* \hookrightarrow {\Bbb C} [M_{p\times r}]$ (resp. $F^{\mu} \hookrightarrow {\Bbb C} [M_{q \times r}]$) associe alors \`a l'\'el\'ement 
$e_T$ (resp. $f_U$) le polyn\^ome
$$\Delta_T (X) := \Delta_{i_{1,1}, \ldots , i_{a,1}} (X) \ldots \Delta_{i_{1, q-b} , \ldots , i_{a, q-b}}  (X) \in {\Bbb C} [M_{p\times r}]$$
$$\left( {\rm resp.} \  
\widetilde{\Delta}_U (Y) := \widetilde{\Delta}_{i_{1,1}, \ldots , i_{b,1}} (Y) \ldots \widetilde{\Delta}_{i_{1,p-a} , \ldots , i_{b,p-a}} (Y) \in {\Bbb C} [M_{q\times r}] \right) ,$$
o\`u 
$$\Delta_{i_1 , \ldots , i_a} (X) = \det \left(
\begin{array}{ccc} 
X_{i_1 , 1} & \cdots & X_{i_1 , a} \\
\vdots & & \vdots \\
X_{i_a , 1} & \cdots & X_{i_a , a} 
\end{array} \right) \in {\Bbb C} [M_{p\times r}] \ \mbox{ pour } 1 \leq i_1 < \ldots < i_a \leq p $$
et 
$$\widetilde{\Delta}_{i_1 , \ldots , i_b} (Y) = \det \left(
\begin{array}{ccc} 
Y_{q-i_b +1 , 1} & \cdots & Y_{q- i_b +1, b} \\
\vdots & & \vdots \\
Y_{q-i_1 +1 , 1} & \cdots & Y_{q -i_1+1 , b} 
\end{array} \right) \in {\Bbb C} [M_{q\times r}] \ \mbox{ pour } 1 \leq i_1 < \ldots < i_b \leq q .$$

\begin{thm} \label{Tc}
La classe 
\begin{eqnarray} \label{classe}
\psi^{(\nu)} \in [{\cal F}^{\infty} \otimes \Omega^{0,aq+ bp -ab} (D ) ]^{SU(p,q)} \cong [{\cal F}^{\infty} \otimes \bigwedge {}^{0, aq+bp-ab} \mathfrak{p}^* ]^{S(U(p) \times U(q))}
\end{eqnarray}
qui \`a un \'el\'ement de la forme $e_T^* \otimes f_U \otimes d \in (E^{\lambda})^* \otimes F^{\mu} \otimes D_r (F)$ 
associe $\Delta_T \cdot \widetilde{\Delta}_U$, d\'efinie une forme diff\'erentielle anti-holomorphe ferm\'ee (harmonique)
non nulle.
\end{thm}
{\it D\'emonstration.} La classe $\psi^{(\nu)}$ appartient \`a l'espace 
$${\rm Hom}_{S(U(p) \times U(q))} (\bigwedge  {}^{0,aq+bp-ab} \mathfrak{p} , {\cal F}^{\infty} ).$$ 
Son image est irr\'eductible; c'est le plus bas $K$-type du sous-$SU(p,q)$-module $\overline{A(\nu)}$ dans ${\cal F}^{\infty}$. Ce module est cohomologique. La forme diff\'erentielle
$$g \cdot \lambda \mapsto \omega (g) \varphi (\lambda )   \  (g\in SU(p,q) , \  \lambda \in \bigwedge {}^{0,aq+bp-ab} \mathfrak{p} )$$
est donc harmonique. Elle est bien s\^ur non nulle. Elle d\'efinie donc bien une forme diff\'erentielle anti-holomorphe ferm\'ee non nulle.

\medskip

Ces formes sont $K$-isotypiques~: elles sont r\'ealis\'ees dans le $K$-type minimal de la repr\'esentation cohomologique, {\it i.e.}
l'image de $\psi$ est $\overline{V(\nu)}$. Lorsque $b=0$ et $a=r=1$, on peut facilement expliciter ces ``classes de Fock''. Supposons donc dor\'enavant $b=0$ et $a=r=1$; alors $\nu = (q)=  1 \times q$.

L'espace tangent anti-holomorphe au point $0$ de l'espace sym\'etrique $D$ s'identifie, comme $U(p) \times U(q)$-module, \`a $E^* \otimes F$. 
Notons $\overline{\xi}_{i,j}$ la forme lin\'eaire 
$$ \overline{\xi}_{i,j} := e_i \otimes f_j^* : \mathfrak{p}^- \cong E^* \otimes F \rightarrow {\Bbb C}$$ 
duale au vecteur $e_i^* \otimes f_j$. (Nous
noterons $\xi_{i,j}$ la forme lin\'eaire $e_i^* \otimes f_j$.) 

\noindent
Les formes $\xi_{ij}$ et $\overline{\xi}_{ij}$ d\'efinissent des formes diff\'erentielles holomorphes et 
anti-holomorphes sur $D$. Remarquons que l'espace $E^{(q)}$ est isomorphe au produit sym\'etrique 
${\rm Sym}^q (E)$; on a donc $\overline{V((q))} = {\rm Sym}^q (E^*) \otimes \bigwedge^q F$. 
Un tableau de Young $T$ de diagramme $(q)$ comme ci-dessus correspond au choix de $q$ entiers 
$1 \leq i_1 \leq \ldots \leq i_q \leq p$; les \'el\'ements 
$$\overline{\xi}_T = \sum_{\sigma \in \mathfrak{S}_{q}} {\rm sgn} (\sigma ) \cdot \overline{\xi}_{i_{1},\sigma (1)}  \wedge \ldots \wedge   \xi_{i_{q}, 
\sigma (q)} $$
forment donc une base duale de $\overline{V((q))}$. Notons enfin $\xi_T$ l'analogue holomorphe de la forme $\overline{\xi}_T$.

La classe de Fock (\ref{classe}) s'\'ecrit alors dans ce cas
$$\psi^{(q)} = \sum_T \Delta_T \overline{\xi}_T,$$
o\`u $T$ parcourt l'ensemble des tableaux de Young de diagramme $(q)$; soit
\begin{eqnarray*}
\psi^{(q)} &  = & \sum_{1 \leq i_1 \leq \ldots \leq i_q \leq p} X_{i_1 , 1} \cdots X_{i_q , 1} \left(  \sum_{\sigma \in \mathfrak{S}_q} {\rm sgn} (\sigma ) \cdot \overline{\xi}_{i_1 ,\sigma (1)} \wedge \ldots \wedge 
\overline{\xi}_{i_q , \sigma (q)} \right) \\
& = & \frac{1}{q!}  \sum_{\sigma , \tau \in  \mathfrak{S}_q} {\rm sgn} (\sigma  \tau^{-1}) \cdot \\
&    &  \sum_{1 \leq i_1 \leq \ldots \leq i_q \leq p} X_{i_{\tau (1)} , 1} \cdots X_{i_{\tau (q)} , 1} 
\overline{\xi}_{i_{\tau (1) } ,\sigma \tau^{-1}  (\tau (1))} \wedge \ldots \wedge \overline{\xi}_{i_{\tau (q)} , \sigma \tau^{-1} (\tau (q))} \\
& = & \frac{1}{q!}  \sum_{1 \leq i_1 , \ldots , i_q \leq p} X_{i_1 , 1} \cdots X_{i_q , 1}  \cdot \left( \sum_{\sigma \in  \mathfrak{S}_q} {\rm sgn} (\sigma ) \cdot 
 \overline{\xi}_{i_1 ,\sigma (1)} \wedge \ldots \wedge \overline{\xi}_{i_q , \sigma (q)} \right) \\
& = & \sum_{1\leq i_1, \ldots , i_q \leq p} X_{i_1 , 1} \cdots X_{i_q , 1} \overline{\xi}_{i_1 ,1} \wedge \ldots \wedge \overline{\xi}_{i_q , q}.
\end{eqnarray*}

\begin{prop} \label{exp}
Lorsque $b=0$ et $a=r=1$, la classe de Fock (\ref{classe}) s'\'ecrit explicitement
$$\psi^{(q)} = \sum_{1\leq i_1, \ldots , i_q \leq p} X_{i_1 , 1} \cdots X_{i_q , 1} \overline{\xi}_{i_1 ,1} \wedge \ldots \wedge \overline{\xi}_{i_q , q}.$$
\end{prop}
 
Au prix de formules peu maniables, on pourrait g\'en\'eraliser cette construction pour $b=0$ et $a=r$ quelconque et obtenir une expression explicite des formes de Fock $\psi^{(q^r)}$. 
On pr\'ef\`ere former le cup-produit  des formes $\psi^{(q)}$, \`a l'aide des groupes en balance
\begin{eqnarray} \label{B1}
\begin{array}{ccc}
U(r) & & U(p,q) \times \ldots \times U(p,q) \\
 | & \times & | \\
U(1)  \times \ldots \times U(1) & & U(p,q) 
\end{array}
\end{eqnarray}
(la deuxi\`eme inclusion verticale est diagonale).

Le diagramme de groupes en balance (\ref{B1}) correspond \`a la d\'ecomposition de l'espace de Fock
$${\cal P} (M_{p\times r} \oplus M_{q \times r}) \cong \underbrace{{\cal P} (M_{p\times 1} \oplus M_{q\times 1}) \otimes \ldots \otimes {\cal P} (M_{p\times 1} \oplus M_{q\times 1})}_{r \ {\rm fois}}.$$

On peut alors prendre pour forme de Fock
$$\psi^{(0,rq)} := \psi_1 \wedge \ldots \wedge \psi_r ,$$
o\`u 
$$\psi_k = \sum_{1\leq i_1, \ldots , i_q \leq p} X_{i_1 , k} \cdots X_{i_q , k} \overline{\xi}_{i_1 ,1} \wedge \ldots \wedge \overline{\xi}_{i_q , q}.$$

\begin{thm} \label{T33}
La  famille de formes de Fock non nulles
$$\psi^{(0,rq)}  \in [{\cal F}^{\infty} \otimes \Omega^{0,rq} (D) ]^{U(p,q)} \cong [{\cal F}^{\infty} \otimes \bigwedge {}^{0,rq} \mathfrak{p}^* ]^{U(p) \times U(q)}$$
($0 \leq r \leq p$) v\'erifie 
\begin{enumerate}
\item $d\psi^{(0,rq)} =0$, i.e. pour tout $(X,Y) \in M_{p\times r} \oplus M_{q\times r}$, $\psi^{(0,rq)} (X,Y)$ est une $rq$-forme anti-holomorphe ferm\'ee sur $D$ qui est 
$G_{(X,Y)}$-invariante;
\item les formes de Fock sont compatibles avec le cup-produit
$$\psi^{(0,r_1q)} \wedge \psi^{(0,r_2q)}  = \psi^{(0,(r_1 + r_2)q)},$$
avec $\psi^{(0,rq)} = 0$ si $r>p$;
\item la classe 
$$[\psi^{(0,rq)} ] \in H^{rq} (\mathfrak{g} , K ; {\cal F}^{\infty} )$$
est non nulle, fortement primitive et $U(p)\times U(q)$-isotypique.
\end{enumerate}
\end{thm}
{\it D\'emonstration.} Cela d\'ecoule des r\'esultats ci-dessus et de \cite[Th\'eor\`eme 5.8]{Lefschetz} pour le point {\it 3.}

\medskip

Remarquons que contrairement aux formes de Fock $\psi^{(q^r)}$ les formes $\psi^{(0,rq)}$ ne sont pas $U(p)\times U(q)$-isotypiques; la forme $\psi^{(q^r)}$ co\"{\i}ncide d'ailleurs avec la partie
$\overline{V((q^r))}$-isotypique de $\psi^{(0,rq)}$. La classe de cohomologie $[\psi^{(0,rq)} ] $ reste isotypique.

Enfin, on aurait pu consid\'erer la polarisation de Fock correspondant au sous-espace anti-holomorphe. On obtient alors l'analogue holomorphe de nos classes de Fock 
en rempla\c{c}ant les polyn\^omes par leurs conjugu\'es complexes et les $\overline{\xi}_{ij}$ par leurs analogues holomorphe. 
Notons ces formes $\psi^{(rq,0)}$.

\subsection{Classes de Fock dans les groupes unitaires et orthogonaux}

Consid\'erons plus g\'en\'eralement la paire r\'eductive duale $(U(r,s), U(p,q))$. Supposons $r=s+k \geq s$. 
On pr\'ef\`ere r\'ealiser le groupe $U(r,s)$ comme groupe
unitaire $G'$ de l'espace hermitien $V' ={\Bbb C}^{r+s}$ muni de la forme 
hermitienne $(\cdot , \cdot )'$ de matrice 
$$\left( 
\begin{array}{ccc}
0 & 0 & -i1_s \\
0 & 1_k & 0 \\
i1_s & 0 & 0 
\end{array} \right).$$ 
Il lui correspond naturellement la polarisation ${\Bbb C}^{r+s} = ({\Bbb C}^s \oplus {\Bbb R}^k) \oplus (i{\Bbb R}^k \oplus {\Bbb C}^r)$ 
associ\'ee \`a la structure complexe
$$J= \left( 
\begin{array}{ccc}
0 & 0  & -1_s \\
0 & i 1_k & 0 \\
1_s & 0 & 0  
\end{array} \right).$$
Nous conservons par ailleurs les notations des sections pr\'ec\'edentes.
 
L'espace $W = V \otimes_{{\Bbb C}} V'$ est naturellement muni de la forme symplectique $B=-2{\rm Im} \left( (,) \otimes \overline{(,)''} \right)$ 
et le couple $(U(p,q) , U(r,s))$ forme une paire r\'eductive duale dans $Sp(W)$. Nous allons consid\'erer le mod\`ele de Fock  ${\cal F}$
associ\'e \`a la structure complexe d\'efinie positive 
$$J_0 = J \otimes I_{p,q} . $$

Dans ce mod\`ele de Fock, l'espace des vecteurs $\widetilde{U((r+s)n)}$-finis s'identifie avec l'anneau des polyn\^omes ($J$-complexes)
sur le complexifi\'e de l'espace vectoriel r\'eel $M_{p\times s} ({\Bbb C}) \oplus M_{q \times s} ({\Bbb C}) \oplus M_{p \times k } ({\Bbb R}) \oplus M_{p\times k} ({\Bbb R})$
que nous identifions, via l'application
$$X\oplus Y \oplus Z \oplus T \mapsto (X \oplus Y , Z \oplus T , \overline{X} \oplus \overline{Y})$$
\`a l'espace vectoriel complexe $(M_{p\times s}  \oplus M_{q \times s} ) \times ( M_{p \times k }  \oplus M_{p\times k} ) \times (M_{p\times s}  \oplus M_{q \times s} )$.
Via ces identifications, l'espace des polyn\^omes $J$-complexes s'obtient comme produit tensoriel
$${\cal P} [ M_{p\times s} \oplus M_{q \times s}] \otimes {\cal P} [M_{p\times k} \oplus M_{q \times k}] \otimes {\cal P} [ \overline{M_{p\times s} \oplus M_{q \times s}}]$$
des espaces 1) de polyn\^omes holomorphes sur $M_{p\times s} \oplus M_{q \times s}$, 2) de polyn\^omes holomorphes sur 
$M_{p \times k} \oplus M_{q \times k}$ et 3) de polyn\^omes anti-holomorphes sur $M_{p\times s} \oplus M_{q \times s}$. 
Cette d\'ecomposition correspond au diagramme de groupes en balance suivant~:
\begin{eqnarray} \label{B3}
\begin{array}{ccc}
U(r,s) & & U(p,q) \times U(p,q) \times U(p,q) \\
 | & \times & | \\
U(s)  \times U(k) \times U(s) & & U(p,q) 
\end{array}
\end{eqnarray}
(la deuxi\`eme inclusion verticale est diagonale). La diagonale ascendante de ce diagramme correspond au produit tensoriel des representations
de Weil associ\'ees aux paires $(U(s) , U(p,q))$ , $(U(r) , U(p,q))$ et $(U(s), U(p,q))$ la derni\`ere \'etant la conjugu\'ee complexe de la premi\`ere.
D'apr\`es la construction du \S 2.1, il leur correspond \`a chacune une classe de Fock; ces classes de Fock sont respectivement anti-holomorphe, 
anti-holomorphe et holomorphe.
Le diagramme (\ref{B3}) permet de consid\'erer leur cup-produit not\'e $\psi^{(sq, rq)}$.

En inversant les r\^oles de $r$ et $s$ et en conjuguant les formes (anti-)holomorphes, on donne plus g\'en\'eralement un sens 
\`a $\psi^{(rq,sq)}$ pour $r$ et $s$ quelconques. On obtient ainsi le th\'eor\`eme suivant.

\begin{thm} \label{TFockU}
La classe de Fock
$$\psi^{(sq,rq)} = \psi_1 \wedge \ldots \wedge \psi_r \wedge \overline{\psi}_1 \wedge \ldots \wedge \overline{\psi}_s \in 
[{\cal F}^{\infty} \otimes \Omega^{sq,rq} (D) ]^{U(p,q)},$$
$(0 \leq r,s \leq p)$ v\'erifie
\begin{enumerate}
\item $d\psi^{(sq,rq)} =0$;
\item les formes de Fock sont compatibles avec le cup-produit
$$\psi^{(s_1 q , r_1 q)} \wedge \psi^{(s_2 q , r_2 q)} = \psi^{((s_1+s_2)q , (r_1 +r_2)q)},$$
avec $\psi^{(sq,rq)}=0$ si $s$ ou $r>p$;
\item la partie fortement primitive de la classe 
$$[\psi^{(sq,rq)} ] \in H^{(r+s)q} (\mathfrak{g}, K ; {\cal F}^{\infty} )$$
est non nulle pour $r +s \leq p$.
\end{enumerate}
\end{thm}
{\it D\'emonstration.} Cela d\'ecoule imm\'ediatement du th\'eor\`eme \ref{T33} et de \cite[Th\'eor\`eme 5.8]{Lefschetz} pour le dernier point.

\medskip

Supposons $r +s \leq p$. 
La classe de Fock $\psi^{(rq,sq)}$ d\'efinit un \'el\'ement de 
$${\rm Hom}_{U(p) \times U(q)} (\bigwedge {}^{(rq,sq)} \mathfrak{p} , {\cal F}^{\infty});$$
son image n'est pas irr\'eductible; la classe de cohomologie $[\psi^{(rq,sq)}]$ n'est m\^eme pas isotypique. Remarquons
n\'eanmoins que cette derni\`ere poss\`ede une unique composante isotypique fortement primitive. 
Le vecteur $v(r\times q) \otimes w( s \times q)^* \in \bigwedge^{(rq,sq)} \mathfrak{p}$ d\'efini par 
$$v( r \times q ) = \bigwedge_{i=1}^r \bigwedge_{j=1}^q e_i \otimes f_j^* \in \bigwedge {}^{rq} (E \otimes F^* ) = \bigwedge {}^{rq} \mathfrak{p}^+$$
et
$$w( s \times q )^* = \bigwedge_{i=1}^s \bigwedge_{j=1}^q e_{p-i+1}^* \otimes f_{q-j+1} \in \bigwedge {}^{sq} (E^* \otimes F ) = \bigwedge {}^{sq} \mathfrak{p}^- ,$$
engendre en effet, sous l'action de $U(p) \times U(q)$, un sous-module irr\'eductible de $\bigwedge^{(rq,sq)} \mathfrak{p}$. Notons $V((q^r), (q^{p-s}))$ ce module; c'est le plus bas $K$-type d'une 
repr\'esentation cohomologique de $U(p,q)$ de degr\'e fortement primitif $(r+s)q$, le seul qui intervient dans la d\'ecomposition en irr\'eductibles du produit tensoriel 
$$V((q^r)) \otimes \overline{V((q^s))} \mbox{ (resp. } \underbrace{V((q^1)) \otimes \ldots \otimes V((q^1))}_{r \; {\rm fois}} \otimes  
\underbrace{\overline{V((q^1))} \otimes \ldots \otimes \overline{V((q^1))}}_{s \; {\rm fois}}).$$ 
Il intervient dans celui-ci avec multiplicit\'e un. Notons $\psi^{(rq,sq)}_+$ la composante de $\psi^{(rq,rq)}$ relative \`a ce $K$-type. Il serait possible
de d\'ecrire explicitement $\psi^{(rq,sq)}_+$ mais au prix de formules peu maniables.

\medskip
\noindent
{\bf Remarque.} La forme $\psi^{(rq,rq)}_+$ est consid\'er\'ee par Tong et Wang dans \cite{TongWang}; elle n'y est pas d\'ecrite explicitement. 
Une autre fa\c{c}on de construire cette forme, plus dans la mani\`ere de Tong et Wang, 
est de consid\'erer le diagramme de groupes en balance suivant~:
\begin{eqnarray} \label{B13}
\begin{array}{ccc}
U(r,r) \times U(r,r) & & U(p,q)  \\
 | & \times & | \\
U(r,r)  & & U(p) \times U(q) 
\end{array}
\end{eqnarray}
(la premi\`ere inclusion verticale est diagonale). Les correspondances $(U(r,r), U(p))$ et $(U(r,r), U(q))$ sont bien comprises, cf. th\'eor\`eme \ref{KV}, le $U(p)\times U(q)$-module
$V((q^r), (q^{p-r}))$ se d\'ecompose en le produit tensoriel du $U(p)$-module $S^q \left( \bigwedge^r E \right) \otimes S^q \left( \bigwedge^r E^* \right)$ par le $U(q)$-module trivial. 
Chacun de ces modules est en correspondance avec un $U(r,r)$-module holomorphe de plus haut poids. Le diagramme de groupes en balance (\ref{B13}) implique alors que le produit
tensoriel de ces deux modules holomorphes de plus haut poids contient un sous-module irr\'eductible en correspondance avec un $U(p,q)$-module qui contient le $K$-type $V((q^r), (q^{p-r}))$.
En \'etudiant cette correspondance Tong et Wang construisent \'egalement (et toujours de mani\`ere non explicite) dans \cite{TongWang} la forme $\psi^{(rq,rq)}_+$.

\bigskip

On passe des groupes unitaires aux groupes orthogonaux en consid\'erant le diagramme de groupes en balance suivant~:
\begin{eqnarray} \label{B5}
\begin{array}{ccc}
Sp(2r) & & U(p,q) \\
 | & \times & | \\
U(r)  & & O(p,q). 
\end{array}
\end{eqnarray}

Dans le mod\`ele de Fock ${\cal P} [ M_{pÊ\times r} ({\Bbb R}) \oplus M_{q \times r} ({\Bbb R})]$  de la repr\'esentation de l'oscillateur, on obtient
imm\'ediatement le th\'eor\`eme suivant o\`u cette fois $D = SO(p,q) /(SO(p) \times SO(q))$ d\'esigne l'espace sym\'etrique (non connexe $\neq D^+$ !) 
associ\'e au groupe orthogonal $O(p,q)$ et ${\cal F}$ le mod\`ele de Fock ci-dessus vu comme repr\'esentation de $\widetilde{Sp(2r)} \cdot O(p,q)$. 

\begin{thm} \label{TFockO}
La classe de Fock
$$\psi^{(rq)} = \psi_1 \wedge \ldots \wedge \psi_r  \in 
[{\cal F}^{\infty} \otimes \Omega^{rq} (D) ]^{O(p,q)},$$
$(0 \leq r \leq p)$ et v\'erifie
\begin{enumerate}
\item $d\psi^{(rq)} =0$;
\item les formes de Fock sont compatibles avec le cup-produit
$$\psi^{(r_1 q)} \wedge \psi^{(r_2 q)} = \psi^{((r_1+r_2)q)},$$
avec $\psi^{(rq)}=0$ si $r>p$;
\item la partie fortement primitive de la classe 
$$[\psi^{(rq)} ] \in H^{rq} (\mathfrak{g}, K ; {\cal F}^{\infty} )$$
est non nulle pour $r \leq p/2$.
\end{enumerate}
\end{thm}
{\it D\'emonstration.} Cela d\'ecoule imm\'ediatement du th\'eor\`eme \ref{T33} et de \cite[Th\'eor\`eme ?]{Lefschetz}.

\medskip

\subsection{Classes de Kudla et Millson}

Dans \cite{KudlaMillson1}, Kudla et Millson construisent des classes de Schwartz assoc\'ees aux paires r\'eductives duales 
$(U(r,r) , U(p,q))$ et $(Sp(2r) , O(p,q))$. Montrons que les classes construites ci-dessus s'identifient aux classes de Kudla et Millson dans un mod\`ele
de Schr\"odinger appropri\'e. Commen\c{c}ons par le cas des groupes unitaires.

En conservant les notations de la section pr\'ec\'edente, on a cette fois $r=s$ et l'on r\'ealise le groupe $U(r,r)$ 
comme groupe unitaire $G'$ de l'espace hermitien $V' ={\Bbb C}^{2r}$ muni de la forme 
hermitienne $(\cdot , \cdot )'$ de matrice 
$$\left( 
\begin{array}{cc}
 & -i1_r \\
i1_r & 
\end{array} \right).$$ Il lui correspond naturellement la polarisation ${\Bbb C}^{2r} = {\Bbb C}^r \oplus {\Bbb C}^r$ associ\'ee \`a la structure complexe
$$J= \left( 
\begin{array}{cc}
  & -1_r \\
1_r & 
\end{array} \right).$$

Remarquons que le sous-groupe compact maximal de $U(r,r)$, isomorphe \`a $U(r) \times U(r)$, est le sous-groupe qui centralise $J$. Il lui
correspond une d\'ecomposition de l'espace $V'= {\Bbb C}^{2r}$ en somme directe de deux sous-espaces complexes de dimension $r$
sur lesquels $J$ induit respectivement la structure complexe usuelle ``multiplication pat $i$'' et l'oppos\'ee de celle-ci.
 
L'espace $W = V \otimes_{{\Bbb C}} V'$ est naturellement muni de la forme symplectique $B=-2{\rm Im} \left( (,) \otimes \overline{(,)''} \right)$ 
et le couple $(U(p,q) , U(r,r))$ forme une paire r\'eductive duale dans $Sp(W)$. Nous allons consid\'erer le mod\`ele de Schr\"odinger 
correspondant \`a la polarisation 
$$(M_{p\times r} \oplus M_{q \times r} ) \oplus  (M_{p\times r} \oplus M_{q \times r} )$$
de $W$, 
associ\'e \`a la structure complexe d\'efinie positive $J_0= J \otimes I_{p,q}$. 
Remarquons que l'espace $M_{p\times r} \oplus M_{q \times r}$ s'identifie naturellement \`a l'espace $V^r$. 

Transposons les r\'esultats de la section pr\'ec\'edente dans le mod\`ele de Schr\"odinger associ\'e.
Le sous-espace ${\rm S} [M_{p\times r} \oplus M_{q \times r}] \cong {\rm S} (V^r )$ de la repr\'esentation de l'oscillateur dans ce mod\`ele de Schr\"odinger est 
\'egale \`a 
$$\varphi^0_+ {\Bbb R} [M_{p\times r} \oplus M_{q \times r}],$$
o\`u $\varphi^0_+$ est la fonction de Schwartz gaussienne de $M_{p\times r} \oplus M_{q \times r}$ dans ${\Bbb R}$ donn\'ee par 
$$\varphi^0_+ (X,Y) = \exp \left( - \pi \sum_{j=1}^r \left( \sum_{i=1}^p |X_{ij}|^2 + \sum_{i=1}^q |Y_{ij}|^2 \right) \right)$$
et ${\Bbb R}[M_{p\times r} \oplus M_{q \times r}]$ d\'esigne l'espace des polyn\^omes sur l'espace vectoriel r\'eel $M_{p\times r} \oplus M_{q \times r}$. 
Alors, si $P({\rm Re} X , {\rm Im} X , {\rm Re} Y , {\rm Im} Y)$ est un polyn\^ome sur l'espace vectoriel r\'eel $M_{p\times r} \oplus M_{q \times r}$, 
il existe un polyn\^ome complexe sur  $(M_{p\times r} \oplus M_{q \times r}) \times (M_{p\times r} \oplus M_{q \times r})$, toujours not\'e $P$, tel que 
$$P(X,Y,\overline{X} , \overline{Y}) = P({\rm Re} X , {\rm Im} X , {\rm Re} Y , {\rm Im} Y).$$

Remarquons maintenant que
$$\psi^{(rq,rq)} = \psi_{1,1} \wedge \ldots \wedge \psi_{r,r}$$
o\`u 
$$\psi_{k,k} = \sum_{
\begin{array}{c}
1 \leq i_1 , \ldots , i_q \leq p \\
1 \leq j_1 , \ldots , j_q \leq p 
\end{array}} (X_1)_{i_1, k } \ldots (X_1)_{i_q , k} \cdot (X_2)_{j_1 , k} \ldots (X_2)_{j_q, k} \cdot
(\overline{\xi}_{i_1, 1} \wedge \xi_{j_1 , 1}) \wedge \ldots \wedge (\overline{\xi}_{i_q , q} \wedge \xi_{j_q , q}) .$$
Notons $D_i$ ($i=1,2$) l'op\'erateur 
$$D^i_{k, j} = \sum_{l=1}^p \left[ \Xi_{l,j}^i \otimes (X_i)_{l,k} \right] ,$$
o\`u $\Xi_{l,j}^i$ d\'esigne l'op\'erateur de multiplication par $\overline{\xi}_{l,j}$, si $i=1$, et par $\xi_{l,j}$, si $i=2$; de telle sorte que 
$$\psi_{k,k} = \prod_{k=1}^r \prod_{j=1}^q D^1_{k,j} D^2_{k,j} \cdot 1.$$

On passe alors du mod\`ele de Fock au mod\`ele de Schr\"odinger en remarquant qu'\`a l'op\'erateur de multiplication par $(X_1)_{i,j}$ (resp. $(X_2)_{i,j}$), dans le mod\`ele 
de Fock, il correspond, dans le mod\`ele de Schr\"odinger (des polyn\^omes r\'eels en $M_{p\times r} \oplus M_{q \times r}$), l'op\'erateur 
$$\frac{\partial }{\partial \overline{X}_{i,j}} - 2\pi X_{i,j}  \  \left( {\rm resp.} \; ( \frac{\partial }{\partial X_{i,j}} -2\pi \overline{X}_{i,j} \right) .$$

On applique alors les calculs de Kudla et Millson. Remarquons que les coordonn\'ees $X_{i,j}$ sont les coordonn\'ees complexes de 
$$M_{(p+q) \times r} \cong V^r.$$
Plus pr\'ecisemment, si $\{ v_1 , \ldots v_m \}$ est une base hermitienne de $V$ telle que 
$$(v_{\alpha } , v_{\alpha} ) =1 , \; 1 \leq \alpha \leq p  \; \mbox{ et } \; (v_{\beta } , v_{\beta} ) =-1 , \; p+1 \leq \beta \leq n+q ,$$
pour $i=1, \ldots , p$ (resp. $j=1 , \ldots , q$) et $k=1, \ldots , r$, $X_{i,k}$ (resp. $Y_{j,k}$) est la $i$-\`eme (resp. $p+j$-\`eme) coordonn\'ee complexe dans le 
$k$-\`eme facteur $V$ de $V^r$, relativement \`a la base $(v_{\alpha})$. 
La gaussienne $\varphi^0_+$ sur $V^r$ co\"{\i}ncide avec celle de Kudla et Millson.

Il correspond aux op\'erateurs $D^1_{k,j}$, les op\'erateurs
$$\bigtriangledown_{k,j} = \frac{1}{2} \sum_{l=1}^p \left[ \overline{\Xi}_{l,j} \otimes  \left(X_{l,k} - \frac{\partial }{\partial \overline{X}_{l,k}}  \right) \right]$$
et aux op\'erateurs $D^2_{k,j}$ les op\'erateurs conjugu\'es $\overline{\bigtriangledown}_{k,j}$. Il correspond donc \`a la classe de Fock $\psi^{(rq,rq)}$ la classe de Schwartz
$$\varphi^{(rq,rq)} = \prod_{k=1}^r \prod_{j=1}^q \bigtriangledown_{k,j} \overline{\bigtriangledown}_{k,j} \cdot \varphi_0.$$
Cette forme est construite et \'etudi\'ee par Kudla et Millson dans \cite{KudlaMillson1}. Ils montrent en particulier la proposition suivante.

\begin{prop}
La forme de Schwartz 
\begin{eqnarray*}
\varphi^{(rq,rq)} & = & \prod_{k=1}^r \prod_{j=1}^q \bigtriangledown_{k,j} \overline{\bigtriangledown}_{k,j} \cdot \varphi_0 \\
& = & \varphi_1 \wedge \ldots \wedge \varphi_r
\end{eqnarray*}
o\`u 
$$\varphi_k = \prod_{j=1}^q  \bigtriangledown_{k,j} \overline{\bigtriangledown}_{k,j} \cdot \varphi_{0,k}$$
et
$$\varphi_{0,k} (X,Y) = \exp \left( - \pi\left(  \sum_{i=1}^p |X_{i,k} |^2  + \sum_{j=1}^q |Y_{j,k} |^2 \right) \right).$$
Et on a 
$$\begin{array}{ccl}
\varphi_k & = & \sum_{\lambda = 0}^q C(q, \lambda ) A[ \omega (k,1) \wedge \overline{\omega} (k,1) \wedge \ldots \wedge \omega (k,q-\lambda ) \wedge \overline{\omega} 
(k,q-\lambda ) \\
& & \wedge \Omega (q-\lambda +1 , q-\lambda +1) \wedge \ldots \wedge \Omega (q,q) ] \cdot \varphi_{0,k} ,
\end{array}$$
o\`u 
$C(q, \lambda ) = \left( \frac{-1}{2\pi} \right)^{\lambda} \frac{(q !)^2}{\lambda ! ((q-\lambda )!)^2 }$, 
$$\omega (k,j) = \sum_{l=1}^p X_{l,k} \overline{\xi}_{l,j}, $$
$$\overline{\omega} (k,j) = \sum_{l=1}^p \overline{X}_{l,k} \xi_{l,j}, $$
$$\Omega (j,j) = \sum_{l=1}^p \overline{\xi}_{l,j} \wedge \xi_{l,j},$$
et 
$$\begin{array}{l}
A[\omega(k,1) \wedge \ldots \wedge \Omega (q,q)] \\
 \   \  = \frac{1}{(q!)^2} \sum_{\sigma , \overline{\sigma} \in \mathfrak{S}_q} {\rm sgn} (\sigma \cdot \overline{\sigma}) \omega (r, \sigma (1)) \wedge 
 \overline{\omega}(r, \overline{\sigma}(1)) \wedge \ldots \wedge \Omega (\sigma (q) , \overline{\sigma} (q)) .
 \end{array}$$
 \end{prop}

Rappelons que le groupe $G' = U(r,r)$ agit sur l'espace de Schwartz ${\cal S} (V^r)$ de $V^r$ via la
repr\'esentation de Weil $\omega$ associ\'ee au caract\`ere additif $x \mapsto e(x) = \exp (2i\pi x)$ de 
${\Bbb R}$. Soit $K' \subset G'$ le sous-groupe compact maximal de $G'$. Le groupe $K' = U(r) \times U(r)$ admet 
un caract\`ere ${\rm det}_+ \otimes {\rm det}_-$. 

\begin{thm}[Kudla-Millson]
L'\'el\'ement $\varphi^{(rq,rq)}$ d\'efinit une classe de Schwartz 
$$\varphi^{(rq,rq)} \in \left[ {\cal S}(V^r) \otimes \Omega^{rq, rq} (D )\right]^{U(p,q)} \cong \left[  {\cal S}(V^r) \otimes \bigwedge {}^{rq, rq} (\mathfrak{p}^* )\right]^{U(p) \times U(q)}$$
telle que 
\begin{enumerate}
\item $d\varphi^{(rq,rq)} =0$, {\it i.e.}, pour tout $x\in V^r$, $\varphi^{(rq,rq)} (x)$ d\'efinit une $(rq, rq)$-forme ferm\'ee sur $D$ qui est $U(p,q)_x$-invariante;
\item sous l'action de $K' = U(r) \times U(r)$ sur ${\cal S} (V^r)$ via la repr\'esentation de Weil,
$\varphi^{(rq,rq)}$ se transforme selon ${\rm det}_+^m \otimes {\rm det}_-^m$;
\item $\varphi^{(rq,rq)} (0) = c(r) \Omega^r$, pour une certaine constante $c(r)$.
\end{enumerate}
\end{thm}

La forme $\varphi^{(rq,rq)}$ d\'efinit donc une classe de $(\mathfrak{g}, K)$-cohomologie dans la repr\'esentation ${\cal S}(V^r )$. Il d\'ecoule du th\'eor\`eme 
\ref{TFockU} que sa composante fortement primitive dans $H^{2rq} (\mathfrak{g}, K ; {\cal S} (V^r))$ est non nulle si et seulement si $r \leq p/2$. La composante
$K$-isotypique de plus haut poids correspond donc \`a la forme de Tong et Wang mentionn\'ee dans la section pr\'ec\'edente.

\medskip

Remarquons que dans le mod\`ele mixte correspondant au diagramme de groupes en balance
\begin{eqnarray} \label{B4}
\begin{array}{ccc}
U(r,s) & & U(p,q) \times U(p,q) \\
 | & \times & | \\
U(s,s)  \times U(k) & & U(p,q), 
\end{array}
\end{eqnarray}
il correspond aux formes $\psi^{(sq,rq)}$ les formes de ``Fock/Schwartz''
\begin{eqnarray} \label{FS}
\begin{array}{cl}
\varphi^{(sq,rq)} & =  \varphi_1 \wedge \ldots \wedge \varphi_s \wedge \psi_{s+1} \wedge \ldots \wedge \psi_{r} \\
 & \in \left[ {\cal S} \left(V^s , {\cal P} [M_{p\times k} \oplus M_{q \times k}]\right) \otimes \Omega^{sq, rq} (X_{p,q} )\right]^G \\
& \cong \left[  {\cal S}\left(V^s , {\cal P} [M_{p\times k} \oplus M_{q \times k}]\right) \otimes \bigwedge {}^{sq, rq} (\mathfrak{p}^* )\right]^{K}
\end{array}
\end{eqnarray}
o\`u l'on identifie ${\rm S}\left(V^s , {\cal P} [M_{p\times k} \oplus M_{q \times k}]\right)$ \`a l'espace  
$$\widehat{\varphi}_0 \cdot {\cal P} \left[ M_{p\times r} \oplus M_{q\times r}\right] ,$$
o\`u
$$\widehat{\varphi}_0 (X,Y) = \exp \left( - \pi \sum_{j=1}^s \left( \sum_{i=1}^p |X_{ij}|^2 + \sum_{i=1}^q |Y_{ij}|^2 \right) \right)$$
et $M_{p\times k} \oplus M_{q \times k}$ est plong\'e dans  
$M_{p\times r} \oplus M_{q\times r}$ via les $k$ derni\`eres lignes. C'est le mod\`ele que nous pr\'ef\`ererons.

\bigskip

Le cas des groupes orthogonaux se traite de la m\^eme mani\`ere. En passant du mod\`ele de Fock au mod\`ele de Schr\"odinger on r\'eobtient les formes de Kudla et Millson.

\begin{prop}
La forme de Schwartz correspondante \`a $\psi^{(rq)}$ est
\begin{eqnarray*}
\varphi^{(rq)} & = & \varphi_1 \wedge \ldots \wedge \varphi_r
\end{eqnarray*}
o\`u 
$$\begin{array}{ccl}
\varphi_k & = & \sum_{\lambda = 0}^{[q/2]} C(q, \lambda ) A[ \omega (k,1) \wedge \wedge \ldots \wedge \omega (k,q-2\lambda )  \\
& & \wedge \Omega (q-2\lambda +1 , q-\lambda +2) \wedge \ldots \wedge \Omega (q-1,q) ] \cdot \varphi_{0,k} ,
\end{array}$$
avec
$$\varphi_{0,k} (X,Y) = \exp \left( - \pi \left( \sum_{i=1}^p |X_{i,k} |^2 + \sum_{j=1}^q |Y_{j,k} |^2 \right) \right),$$
$$C(q, \lambda ) = \left( \frac{-1}{4\pi} \right)^{\lambda} \frac{q !}{2^{\lambda} \lambda ! (q-2\lambda )! },$$
$$\Omega (j-1,j) = \sum_{l=1}^p \xi_{l,j-1} \wedge \xi_{l,j}$$
et 
$$\begin{array}{l}
A[\omega(k,1) \wedge \ldots \wedge \Omega (q-1,q)] \\
 \   \  = \frac{1}{q!} \sum_{\sigma \in \mathfrak{S}_q} {\rm sgn} (\sigma ) \omega (k, \sigma (1)) \wedge 
 \ldots \wedge \omega (k, \sigma (q-2\lambda) ) \\
 \; \; \; \wedge \Omega (\sigma (q-2\lambda+1 , \sigma (q-2\lambda +2)) \wedge \ldots \wedge \Omega (\sigma (q-1) , \sigma (q)) .
\end{array}$$
\end{prop}

Le groupe $G' = Sp(2r)$ agit sur l'espace de Schwartz ${\cal S} (V^r)$ de $V^r$ via la
repr\'esentation de Weil $\omega$ associ\'ee au caract\`ere additif $x \mapsto e(x) = \exp (2i\pi x)$ de 
${\Bbb R}$. Soit $K' \subset G'$ le sous-groupe compact maximal de $G'$. La pr\'eimage du caract\`ere ${\rm  det}$ du groupe $K' = U(r)$ admet
une racine carr\'ee $({\rm det} )^{1/2}$ sur le rev\^etement m\'etaplectique $\widetilde{U}(r)$. Dans le th\'eor\`eme suivant $D=SO(p,q)/SO(p) \times SO(q)$.

\begin{thm}[Kudla-Millson]
L'\'el\'ement $\varphi^{(rq)}$ d\'efinit une classe de Schwartz 
$$\varphi^{(rq)} \in \left[ {\cal S}(V^r) \otimes \Omega^{rq} (D )\right]^{O(p,q)} \cong \left[  {\cal S}(V^r) \otimes \bigwedge {}^{rq} (\mathfrak{p}^* )\right]^{O(p) \times O(q)}$$
telle que 
\begin{enumerate}
\item $d\varphi^{rq)} =0$, {\it i.e.}, pour tout $x\in V^r$, $\varphi^{rq} (x)$ d\'efinit une $rq$-forme ferm\'ee sur $D$ qui est $O(p,q)_x$-invariante;
\item sous l'action de $K' = U(r)$ sur ${\cal S} (V^r)$ via la repr\'esentation de Weil,
$\varphi^{rq}$ se transforme selon ${\rm det}^{m/2}$;
\item $\varphi^{rq} (0) = 0$ si $q$ est impair et 
$$\varphi^{rq} (0) = c(r) \sum_{\sigma \in \mathfrak{S}_q} {\rm sgn} (\sigma ) \Omega_{\sigma (1) , \sigma (2)} \wedge \ldots \wedge \Omega_{\sigma (q-1), \sigma (q)},$$
pour une certaine constante $c(r)$, si $q$ est pair.
\end{enumerate}
\end{thm}

La forme $\varphi^{rq}$ d\'efinit donc une classe de $(\mathfrak{g}, K)$-cohomologie dans la repr\'esentation ${\cal S}(V^r )$; l\`a encore,
sa composante fortement primitive dans $H^{rq} (\mathfrak{g}, K ; {\cal S} (V^r))$ est non nulle si et seulement si $r \leq p/2$.

\bigskip

Les formes ferm\'ees $\varphi^{(rq,sq)}$, que ce soit dans le cas des groupes unitaires ou orthogonaux, v\'erifient donc le point {\it 3.} du th\'eor\`eme \ref{T1}; 
elles v\'erifient le point {\it 1.} par construction. Pour conclure remarquons que si $U \subset V$ est un sous-espace de dimension
$l \leq p$ tel que la restriction de $(.)$ \`a $U$ soit totalement d\'efinie positive et $G_U \subset G$ comme dans l'introduction mais maintenant sur ${\Bbb R}$,
le lemme suivant d\'ecoule de la consid\'eration du diagramme de groupe en balance~:
$$
\begin{array}{ccc}
G' \times G' & & G \\
| & \times & | \\
G' & & G_U \times G_{U^{\perp}}
\end{array}
$$
(la premi\`ere inclusion verticale est diagonale). Ici si $G= U(p,q)$ (resp. $O(p,q)$), $G_U = U(p-l,q)$ (resp. $O(p-l,q)$) et $G_{U^{\perp}} = U(l)$ (resp. $O(l)$) et 
$G' = U(r,s)$ (resp. $Sp(2r)$).

\begin{lem} \label{restri}
Soit $\varphi_+^0$ la gaussienne associ\'ee \`a la paire duale $(G' , G_U)$, comme au-dessus. Alors l'application de restriction
${\rm res} : \Omega^* (D) \rightarrow \Omega^* (D_U )$, sur les formes diff\'erentielles, envoie $\varphi^{(rq,sq)}$ sur 
$$\varphi_+^0 \otimes \varphi_{U^{\perp}}^{(rq,sq)},$$
o\`u $\varphi_{U^{\perp}}^{(rq,sq)}$ est la forme construite ci-dessus pour $U^{\perp}$ et $G_U$.
\end{lem}

Le lemme \ref{restri} pr\'ecise le point {\it 2.} du th\'eor\`eme \ref{T1} et en conclut la d\'emonstration.

\section{Cup-produits de s\'eries th\^eta}

Revenons maintenant \`a la situation globale de l'introduction. 
Rappelons que $V_K '$ est un espace vectoriel de dimension fini sur $K$, $(,)'$ une forme symplectique (resp. anti-hermitienne) non d\'eg\'en\'er\'ee sur $V_K '$
et $G'$ le groupe r\'eductif sur ${\Bbb Q}$ obtenu, par 
restriction des scalaires de $K$ (resp. $K_0$) \`a ${\Bbb Q}$, \`a partir du groupe des isom\'etries de $(.)'$ et que l'on suppose que $G^{' (1)} \cong Sp(2r)$ (resp.
$G^{' (1)} \cong U(r,s)$) et $V_K '$ est d'indice de Witt maximal (\'egal \`a $|r-s|$) sous cette condition. \'Etant donn\'e que $r$ et $s$ seront amen\'es \`a varier,
nous noterons $G'_{(r,s)}$ le groupe ci-dessus. 

Soient $(r_i , s_i )$, $i=1,2$, deux couples d'entiers naturels. L'homomorphisme ${\i}_0$ usuel $Sp(2r_1) \times Sp(2r_2 ) \rightarrow Sp(2(r_1 + r_2))$ (resp.
$U(r_1 , s_1 ) \times U(r_2 , s_2 ) \rightarrow U(r_1 + r_2 , s_1 + s_2 )$) se rel\`eve en un homomorphisme
\begin{eqnarray} \label{mor}
\begin{array}{ccc}
\widetilde{G}_{(r_1 , s_1)} ' ({\Bbb A}) \times \widetilde{G}_{(r_2 , s_2)} ' ({\Bbb A}) & \stackrel{\tilde{{\i}}_0}{\rightarrow} & \widetilde{G}_{(r_1+r_2 , s_1 +s_2 )} ' ({\Bbb A}) \\
\downarrow & & \downarrow \\
G_{(r_1 , s_1)} ' ({\Bbb A}) \times G_{(r_2 , s_2)} ' ({\Bbb A}) & \stackrel{\tilde{{\i}}_0}{\rightarrow} & G_{(r_1 +r_2 , s_1 +s_2 )} ' ({\Bbb A}) .
\end{array}
\end{eqnarray}

Fixons un caract\`ere non trivial de ${\Bbb A}_{K_0} /K_0$ et notons ${\cal S} (X_i ({\Bbb A}))$, $i=1,2$, l'espace des vecteurs lisses d'une r\'ealisation (dans un mod\`ele
de Schr\"odinger) de la repr\'esentation de Weil de $\widetilde{G}_{(r_i , s_i )} ' ({\Bbb A})  \cdot G ({\Bbb A})$. Rappelons que l'on a associ\'e \`a une fonction
$\varphi_i \in {\cal S} (X_i ({\Bbb A}_f) )$ la forme diff\'erentielle ferm\'ee $\theta (g_i ' , \varphi_i )$ ($g_i ' \in  \widetilde{G}_{(r_i , s_i )} ' ({\Bbb R})$).
Posons $\varphi = \varphi_1 \otimes \varphi_2 \in {\cal S} ((X_1 \oplus X_2 ) ({\Bbb A}_f ))$. Il d\'ecoule alors du point 1. du th\'eor\`eme \ref{T1} que, 
pour $g_i ' \in  \widetilde{G}_{(r_i , s_i )} ' ({\Bbb R})$, $i=1,2$, on a
\begin{eqnarray} \label{cp}
\theta (\tilde{{\i}}_0 (g_1 ' , g_2 ' ) , \varphi ) = \theta (g_1 ' , \varphi_1 ) \wedge \theta (g_2 ' , \varphi_2 ).
\end{eqnarray}

Le but de cette section est la d\'emonstration des th\'eor\`emes \ref{T2} et \ref{T3}. 
Nous allons en fait montrer que si 
$$Sh(G) := \lim_{\leftarrow_{K}} Sh(G)_K = G({\Bbb Q}) \backslash (D \times G({\Bbb A}_f ))$$
et $r=p$ et $s=0$ (resp. $s=p$) la classe (de degr\'e maximal)
$$[\theta (g' , \varphi )] \in H^{(r+s)q} (Sh (G)) =  \lim_{\leftarrow_{K}} H^{(r+s)q} (Sh(G)_K )$$
est non nulle.  

Quitte \`a reprendre des cup-produits \`a l'aide de (\ref{cp}) et \`a translater par un \'el\'ement de $G({\Bbb A}_f)$ pour intervenir dans la cohomologie de 
$Sh^0 (G)$, le th\'eor\`eme \ref{T2}. Le th\'eor\`eme s'en d\'eduit pareillement en utilisant le point {\it 2.} du th\'eor\`eme \ref{T1}.

Dans la suite nous supposerons
donc $G' = G'_{(r,s)}$ avec $r=p$ et $s=0$ (resp. $s=p$). Les formes que nous consid\'erons sont alors celles consid\'er\'ees par Kudla et Millson. La notation
$[\theta (g' , \varphi )]$ d\'esigne dor\'enavant une classe dans $H^* (Sh(G))$.

\subsection{Quelques r\'esultats de Kudla et Millson}

\paragraph{Lemme de Thom.} Soit $r$ un entier naturel $\leq p$. 
Notons $\varphi^{(r)}$ la forme de Kudla-Millson $\varphi^{(rq,rq)}$ dans le cas unitaire $G^{(1)}=U(p,q)$ et 
$\varphi^{(rq)}$ dans le cas orthogonal $G^{(1)}=O(p,q)$. Dans ce paragraphe les facteurs compacts de $G$ ne jouent aucun r\^ole, on note
$$G = G^{(1)} \times U$$
et 
$$V = V^{(1)} \oplus Z,$$
o\`u $U$ est compact et $(,)_{|Z}$ est d\'efinie positive. Soit $\varphi \in {\cal S} (Z^r)$ la fonction
$$\varphi ' (z_1 , \ldots , z_r ) = e^{-\pi \sum_{i=1}^r (z_i , z_i )}.$$
La multiplication par $\varphi '$ identifie $H^* (\mathfrak{g} , K ; {\rm S} ((V^{(1)})^r))$ et $H^* (\mathfrak{g} , K ; {\rm S} (V^r))$. On peut donc
supposer $\mu =1$ et travailler directement avec $\varphi^{(r)} \in H^* (\mathfrak{g} , K ; {\rm S}(V^r ))$.
Remarquons que $\varphi^{(r)}$ est multipli\'ee, sous l'action de $\widetilde{U} (r)\subset \widetilde{Sp}(2r)$ (resp. $U(r) \times U(r) \subset U(r,r)$), par 
le caract\`ere $(\det )^{(p+q)/2}$ (resp. $\det_+^{p+q} \otimes \det_-^{-p-q}$). 

Soit $V_+$, $(,)_+$ un espace orthogonal (resp. hermitien) positif de dimension $n=p+q$ sur ${\Bbb R}$ (resp. ${\Bbb C}$) et soit $\varphi_+^0 \in 
{\rm S} (V_+^r )$ la gaussienne 
$$\varphi_+^0 (x) = \exp \left( -\pi \sum_{i=1}^{n} |x_i |^2 \right).$$
Sous l'action de la repr\'esentation de Weil $\omega_+$ de $\widetilde{Sp}(2r)$ (resp. $U(r,r)$) associ\'ee \`a $V_+$, l'action d'un \'el\'ement $k'$ du compact
maximal $K'$ est alors donn\'ee par~:
\begin{eqnarray} \label{actK}
\omega_+ (k ' ) \varphi_+^0 = \det (k' )^{(p+q)/2} \varphi_+^0  \;  \left( {\rm resp.} \; = \det (k'_+ )^{p+q} \otimes \det (k'_- )^{-p-q} \varphi_-^0 \right).
\end{eqnarray}
De plus, si $x \in V_+^r$ avec $(x,x)/2 = \beta$ matrice symm\'etrique (resp. hermitienne) de taille $r\times r$, alors pour $g' \in \widetilde{Sp}(2r)$ (resp. $U(r,r)$),
la fonction (de Whittaker g\'en\'eralis\'ee)
\begin{eqnarray} \label{W}
W_{\beta} (g' ) = \omega_+ (g' ) \varphi_+^0 (x)
\end{eqnarray}
est explicitement calcul\'ee dans \cite[\S 6.6]{KudlaMillson2}. Posons $\tau = u+iv = g' (i \cdot 1_r ) $; c'est un \'el\'ement de l'espace de Siegel de genre $r$ (resp.
de l'espace sym\'etrique associ\'e \`a $U(r,r)$). Quitte \`a multiplier $g'$ par un certain $k'$ \`a droite, on peut
supposer $g' = n' (b) m' (a) $ (d\'ecomposition d'Iwasawa) 
avec $a \in GL(r)$, et $ \det (a) >0$ dans le cas orthogonal, on obtient alors~:
\begin{eqnarray} \label{W2}
W_{\beta } (g' ) = |\det a|^{(p+q)/2} \cdot \exp \left(  {\rm tr} \beta \tau \right) .
\end{eqnarray} 
L'important pour la suite est que cette fonction ne d\'epend que de $\beta$ et $g'$.
%ces calculs sont a mettre en relation avec les formules de la section 2

On d\'efinit maintenant la classe d'Euler (resp. $q$-i\`eme classe de Chern) $c_q$ sur $D^+$ par la formule 
$$c_q = \left\{
\begin{array}{ll}
0 & \mbox{si } q \mbox{ est impair} \\
\frac{1}{l!} \sum_{\sigma \in \mathfrak{S}_q} {\rm sgn} (\sigma ) \Omega (\sigma (1) , \sigma (2)) \wedge \ldots \wedge \Omega (\sigma (2l-1) , \sigma (2l)) & \mbox{si } q=2l
\end{array} \right. $$
$$\left( {\rm resp.} \; 
c_q = \frac{1}{q!} \sum_{\sigma , \overline{\sigma} \in \mathfrak{S}_q} {\rm sgn} (\sigma \cdot \overline{\sigma}) \Omega(\sigma (1), \overline{\sigma}(1)) \wedge \ldots \wedge
\Omega (\sigma(q) , \overline{\sigma} (q)) \right) .$$

Le th\'eor\`eme
suivant est un cas particulier de \cite[Theorem 4.1]{KudlaMillson2} (cf. \'egalement \cite[\S 9]{KudlaMillson3}).

\'Etant donn\'e un \'el\'ement $x \in V^r$, notons $U=U(x)$ le sous-espace engendr\'e par les diff\'erentes composantes de $x$, $G_U=G_x \subset G$ son stabilisateur
dans $G$ et $D_U^+ \subset D^+$ le sous-espace sym\'etrique correspondant. Remarquons que pour tout $g' \in G'$, 
$\omega(g') \varphi^{(r)} (x)$ est une $(rq,rq)$-forme (resp. $rq$-forme) diff\'erentielle ferm\'ee et $G_U$-invariante sur $D$.

\begin{thm} \label{thom}
Soit $\Gamma_U \subset G_U^{{\rm ad} , +}$ un sous-groupe discret cocompact. Alors, pour toute $((p-r)q,(p-r)q)$-forme (resp. $(p-r)q$-forme) $\eta$ ferm\'ee et born\'ee sur 
$\Gamma_U \backslash D^+$, on a 
$$\int_{\Gamma_U \backslash D^+} (\omega (g') \varphi^{(r)} ) (x) \wedge \eta = \varepsilon (x) W_{\beta} (g') \int_{\Gamma_U \backslash D_U^+} c_q^{r-t} \wedge \eta ,$$
o\`u $\varepsilon (x) = \pm 1$ est un signe, toujours \'egal \`a $1$ dans le cas unitaire ou lorsque $q$ est pair et en g\'en\'eral, d\'etermin\'e par $x$,
lorsque $\dim U = r$, $\varepsilon (x) =\pm 1$ selon que $x$ forme une base orient\'ee ou non de $U$, $\beta = \frac12 (x,x)$, 
$W_{\beta} (g') = \omega (g ') \varphi_0 (x)$, $c_q$ est la forme d'Euler (resp. de Chern) et
$t={\rm rang} (\beta) = \dim U$.
\end{thm} 

\paragraph{Cycles sp\'eciaux.} Revenons \`a la situation globale. Soient $\beta \in M_r (K)$ une matrice symm\'etrique (resp. hermitienne) 
que nous supposerons dor\'enavant totalement d\'efinie positive et $\varphi \in {\rm S} (V({\Bbb A}_f )^r)$
une fonction de Schwartz. La fonction $\varphi$ est $K$-invariante pour un certain sous-groupe compact-ouvert $K \subset G({\Bbb A}_f )$.
Dans le cas orthogonal, ni le groupe $G({\Bbb R})$ ni son groupe adjoint $G^{{\rm ad}}({\Bbb R})$ ne sont connexes. Notons 
$G^{{\rm ad}} ({\Bbb R})^+$ la composante connexe de l'identit\'e dans $G^{{\rm ad}} ({\Bbb R})$ et $G({\Bbb R})_+$ sa pr\'eimage dans $G({\Bbb R})$. 
On \'ecrit $G({\Bbb A}_f )= \cup_j G({\Bbb Q})_+ g_j K$, o\`u la r\'eunion (finie) est une union disjointe et $G({\Bbb Q})_+ = G({\Bbb Q}) \cap G({\Bbb R})_+$.
Alors, 
$$Sh (G)_K := G({\Bbb Q}) \backslash (D \times G({\Bbb A}_f )) / K \cong \cup_j \Gamma_j \backslash D^+,$$
o\`u la r\'eunion est disjointe et $\Gamma_j$ est l'image dans $G^{{\rm ad}} ({\Bbb R})^+$ du sous-groupe $\Gamma_j ' = g_j K g_j^{-1} \cap G({\Bbb Q})_+$ de
$G({\Bbb Q})_+$.

Soit $U_K \subset V_K$ un sous-espace rationnel de dimension $r$ et tel que la restriction de $(,)$ \`a $U_K$ soit totalement d\'efinie positive.
Remarquons alors que $r \leq p$.  Notons $U\subset V$ (resp. $U^{\perp} \subset V$) l'espace des points r\'eels 
du ${\Bbb Q}$-espace vectoriel obtenu, \`a partir de $U_K$ (resp. $V_K$), par restriction des scalaires. 
Et $H$ le groupe isomorphe au groupe des points r\'eels du groupe r\'eductif sur ${\Bbb Q}$ obtenu,
par restriction des scalaires de $K$ \`a ${\Bbb Q}$ (resp. $K_0$ \`a ${\Bbb Q}$), \`a partir du groupe des isom\'etries de la forme $(,)_{|U^{\perp}_K}$ sur $U_K^{\perp}$,
de telle mani\`ere qu'on ait un morphisme naturel $H \rightarrow G$. 
Rappelons que l'espace sym\'etrique $D$ s'identifie au sous-ensemble ouvert de la grassmannienne ${\rm Gr}_q (V)$ constitu\'e 
de ceux des $q$-plans $Z$ dans $V$ qui sont tels que $(,)_{|Z}$ est d\'efinie n\'egative. On d\'efinit alors un sous-ensemble $D_U = D_H \subset D$ par 
$$D_U = \{ Z \in D \; : \; Z \subset U^{\perp} \}.$$
Le groupe $H$ est isomorphe au stabilisateur $G_U$ de $U$ dans $G$.

On associe \`a $U$ un {\it cycle connexe} de la mani\`ere suivante. Posons d'abord 
$$\Gamma_{j,U} ' = H({\Bbb Q})_+ \cap g_j K g_j^{-1} = H({\Bbb Q}) \cap \Gamma_j ' .$$
Le groupe $\Gamma_{j,U}'$ est \'egal au stabilisateur de $U$ dans $\Gamma_j '$; notons $\Gamma_{j,U}$ son image dans 
$H^{{\rm ad}} ({\Bbb R})^+$. On a alors une application naturelle
$$\Gamma_{j,U} \backslash D_U \rightarrow \Gamma_j \backslash D.$$
L'image de cette application est un cycle connexe dans $G({\Bbb Q}) \backslash (D \times G({\Bbb A}_f )) / K$ que nous notons
$c(U, g_j , K)$. 

On associe \`a $\beta$ et $\varphi$ la combinaison lin\'eaire suivante de cycles connexes~:
$$Z(\beta , \varphi , K ) := \sum_j \sum_{x \in \Omega_{\beta}  \; {\rm mod} \; \Gamma_j '} \varphi (g_j^{-1} x ) \varepsilon (x) c(U(x) , g_j , K) ,$$
o\`u 
$$\Omega_{\beta} = \left\{ x \in (V_K)^n \; : \; \frac12 (x,x) = \beta \right\}$$
et 
$U(x)$ est le $K$-sous-espace de $V$ engendr\'e par les composantes de $x$. Remarquons que puisque $\beta$ est totalement d\'efinie positive, la restriction
de $(.)$ \`a $U$ est d\'efinie positive. 

Le cycle $Z(\beta , \varphi , K)$ d\'efinit une classe de cohomologie
$$[\beta , \varphi ]^0 := [Z (\beta , \varphi , K) ] \in H^{bq} (Sh (G)_K ) \; \left({\rm resp.} \; H^{(bq,bq)} (Sh (G)_K ) \right),$$ 
o\`u $b= {\rm rang} (\beta)$.

Le cup-produit par la classe d'Euler (resp. $q$-i\`eme classe de Chern) de $D$ induit un op\'erateur
$$H^* (Sh(G)_K ) \rightarrow H^{*} (Sh(G)_K)$$
sur la cohomologie qui commute \`a l'action de $G({\Bbb A}_f )$. On pose alors 
\begin{eqnarray} \label{cd}
[\beta , \varphi] := c_q^{r- {\rm rang} (\beta )} \cdot [\beta , \varphi ]^0 \in H^{rq} (Sh (G)_K ) \; \left({\rm resp.} \; H^{(rq,rq)} (Sh (G)_K ) \right).
\end{eqnarray}

\paragraph{Formes duales.} On peut maintenant faire le lien entre les classes $[\theta (g', \varphi)] \in H^* (Sh(G))$ et les cycles construits ci-dessus.

\'Etant donn\'es $g' \in \widetilde{G} '({\Bbb R})$ et $\beta \in M_r (K)$ sym\'etrique (resp. hermitienne) totalement positive, $\beta \geq 0$, posons
$$W_{\beta} (g ') = W_{\beta^{\sigma_1}} (g_1 ') \ldots W_{\beta^{\sigma_{\mu}}} (g_{\mu} ') .$$

\begin{thm} \label{FD}
\'Etant donn\'es $g' \in \widetilde{G} '({\Bbb R}) \subset \widetilde{G}' ({\Bbb A})$ et $\varphi \in {\rm S} (V({\Bbb A}_f )^r )$, 
$$[\theta ( g' , \varphi ) ] = \sum_{\beta \geq 0} [\beta , \varphi ] W_{\beta} (g') .$$
\end{thm}
{\it D\'emonstration.} Soit $K \subset G({\Bbb A}_f )$ un sous-groupe compact ouvert tel que $\varphi$ soit $K$-invariante et tel que 
l'image de $\Gamma ' = K \cap G({\Bbb Q})_+$ dans $G^{{\rm ad}} ({\Bbb R})$ pr\'eserve l'orientation. Soit $\eta$ une forme ferm\'ee sur $Sh(G)_K$ de degr\'e
$r-p$ (resp. de bidegr\'e $(r-p,r-p)$). Remarquons que 
$$\Omega^* (Sh(G)_K ) \cong \left[ \Omega^* (D) \otimes C^{\infty} (G({\Bbb A}_f ))\right]^{G({\Bbb Q}) \times K} \cong \bigoplus_j \Omega^* (D^+ )^{\Gamma_j},$$
o\`u le second isomorphisme est induit par l'\'evaluation en $g_j$.
Alors,
\begin{eqnarray} \label{86} 
\begin{array}{rl}
\int_{Sh(G)_K} \theta (g' , \varphi ) \wedge \eta & \\
 = & \sum_j \int_{\Gamma_j \backslash D^+} \theta (g', g_j , \varphi ) \wedge \eta (g_j ) \\
= & \sum_j \sum_{\beta } \sum_{x \in \Omega_{\beta} (K) \; {\rm mod} \; \Gamma_j '} \int_{\Gamma_{j,x} \backslash D^+} \omega (g' ) \widetilde{\varphi} (g_j^{-1} x) \wedge \eta (g_j).
\end{array}
\end{eqnarray}
Ici on a utilis\'e le fait que, si $\Gamma_{j,x}$ est l'image dans $\Gamma_j$ du stabilisateur $\Gamma_{j,x} '$ de $x$ dans $\Gamma_j '$, alors 
$$\Gamma_{j,x} ' \backslash \Gamma_j ' \cong \Gamma_{j,x} \backslash \Gamma_j .$$
Le groupe $\Gamma_{j,x} '$ contient en effet le noyau de la projection de $\Gamma_j '$ dans $G^{{\rm ad}} ({\Bbb R})$ et 
l'image de $\Gamma_j '$ dans $G^{{\rm ad}} ({\Bbb R})$ pr\'eserve l'orientation et co\"{\i}ncide donc avec sa projection $\Gamma_j$ dans
$G^{{\rm ad}} ({\Bbb R})^+$.

Le r\'esultat principal de \cite{KudlaMillson3} est que les termes de (\ref{86}) tels que $\beta$ ne sont pas positifs (semi-d\'efinis) sont nuls. Le 
th\'eor\`eme \ref{thom} implique alors que (\ref{86}) est \'egal \`a 
\begin{eqnarray} \label{87}
\begin{array}{l}
\sum_{\beta \geq 0} \sum_j \sum_{x \in \Omega_{\beta} (K) \; {\rm mod} \; \Gamma_j '} \varphi (g_j^{-1} x) \cdot \int_{\Gamma_{j,x} \backslash D^+}
\omega (g') \varphi^{(r)} (x) \wedge \eta (g_j) \\
= \sum_{\beta \geq 0} \sum_j \sum_x \varphi (g_j^{-1} x) \varepsilon (x) \cdot \int_{\Gamma_{j,x} \backslash D_{U(x)}^+} c_q^{r-t} \wedge \eta (g_j) \cdot
W_{\beta} (g') , 
\end{array}
\end{eqnarray}
o\`u $t$ est le rang de $\beta$. Par d\'efinition de $[\beta , \varphi]$, on obtient
$$\int_{Sh(G)_K} \theta (g' , \varphi ) \wedge \eta = \sum_{\beta \geq 0 } \langle [\beta , \varphi ] , \eta \rangle \cdot W_{\beta} (g' ),$$
comme annonc\'e.

\medskip

\subsection {Formule de Siegel-Weil}

Nous rappelons maintenant la formule de Siegel-Weil telle qu'\'etendue par Kudla et Rallis \cite{KudlaRallis} puis Sweet \cite{Sweet}
et Ichino \cite{Ichino}. 

Conservons les notations pr\'ec\'edentes et supposons $r=p$ (le groupe $G^{' (1)}$ est donc isomorphe au groupe $U(p,p)$).
Pour toute fonction de Schwartz $f \in {\rm S} (V({\Bbb A})^p)$, et pour $g' \in \widetilde{G} ' ({\Bbb A})$ et $g  \in G({\Bbb A})$, notons
$$\theta (g', g ; f ) = \sum_{x \in V (K)^p} \omega (g') f(g^{-1} x)$$
la fonction th\^eta usuelle. Soient $P$ le parabolique de Siegel de $G'$, $\rho = (p+1)/2$ (resp. $\rho = p/2$) et $s_0 = (p+q)/2 - \rho$. 
Pour $s \in {\Bbb C}$ soit
$$\Phi (g',s) = \omega (g') f(0) \cdot |a(g')|^{s-s_0},$$
o\`u nous renvoyons aux articles \cite{KudlaRallis}, \cite{Sweet} et \cite{Ichino} pour la d\'efinition de $|a(g')|$. Alors la s\'erie d'Eisenstein, associ\'ee
\`a $f$,
$$E(g',s;f) = \sum_{\gamma \in P({\Bbb Q}) \backslash G' ({\Bbb Q})} \Phi (\gamma g' ,s)$$
est absolument convergente dans le demi-plan Re$(s) > \rho$ et la fonction $s \mapsto E(g',s;f)$ admet un prolongement m\'eromorphe \`a tout le
plan des $s\in {\Bbb C}$. 

Le th\'eor\`eme suivant est d\'emontr\'e dans \cite{KudlaRallis} dans le cas orthogonal et $p+q$ pair, dans \cite{Sweet} dans le cas orthogonal et
$p+q$ impair et dans \cite{Ichino} dans le cas unitaire. 

\begin{thm} \label{SW}
La fonction $s \mapsto E(g',s;f)$ est holomorphe en $s=s_0$ et 
$$E(g',s_0 , f) = \int_{U(V)(K) \backslash U(V)({\Bbb A}_K )} \theta (g', g ; f) dg,$$
o\`u $dg$ est la mesure invariante normalis\'ee de telle mani\`ere que 
$${\rm vol} (U(V)(K) \backslash U(V)({\Bbb A}_K ))=1.$$
\end{thm}

Pour appliquer ce th\'eor\`eme \`a nos s\'eries th\^eta, il nous faut relier l'int\'egrale sur $U(V)(K) \backslash U(V)({\Bbb A}_K )$ \`a l'int\'egrale sur 
$Sh(G)_K$. Pour tout compact ouvert $K$, il existe une application naturelle $H^{pq} (Sh(G)_K ) \rightarrow {\Bbb C}$ (resp. 
$H^{(pq,pq)} (Sh(G)_K) \rightarrow {\Bbb C}$) induite, au niveau des formes diff\'erentielles et apr\`es le choix d'une orientation de $D^+$, 
par l'int\'egration le long de $Sh(G)_K$. Si $K' \subset K$ est un autre compact ouvert, on a alors 
\begin{eqnarray} \label{101}
\begin{array}{lcl}
H^* (Sh(G)_{K'} ) & \rightarrow & {\Bbb C} \\
{\rm pr}^* \uparrow & & \uparrow {\rm deg (pr)} \\
H^* (Sh(G)_K ) & \rightarrow & {\Bbb C} ,
\end{array}
\end{eqnarray}
o\`u deg(pr) est le degr\'e du rev\^etement $Sh(G)_{K'} \rightarrow Sh(G)_K$.

\begin{lem} \label{L101}
On a deg(pr)$=|K/K'(K\cap Z({\Bbb Q}))|$, o\`u $Z$ est le noyau (contenu dans le centre) de la projection $G \rightarrow G^{{\rm ad}}$. 
\end{lem}
{\it D\'emonstration.} La pr\'eimage par pr d'un point $G({\Bbb Q}) (z,g) K \in Sh(G)_K$ est constitu\'e es doubles classes $G({\Bbb Q}) (z,g) k K'$ o\`u $k$
parcourt $K$. Et, 
\begin{eqnarray*}
G({\Bbb Q}) (z,g) k_1 K' = G({\Bbb Q}) (z,g) k_2 K' 
\end{eqnarray*}
si et seulement si $\gamma z = z$ et $\gamma g k_2 k' = g k_1$, pour certains $\gamma \in G({\Bbb Q})$ et $k' \in K'$. On a alors n\'ecessairement 
$\gamma \in G({\Bbb Q})_+$ et $\gamma \in gKg^{-1}$ autrement dit
$$\gamma \in gKg^{-1} \cap G({\Bbb Q})_+ = \Gamma_g ' .$$
On peut supposer $\Gamma_g$, l'image de $\Gamma_g '$ dans $G({\Bbb Q}) / Z({\Bbb Q}) \subset G^{{\rm ad}} ({\Bbb R})$,  
sans torsion (ou consid\'erer un \'el\'ement $z$ ``g\'en\'erique''); son action sur $D$ est alors sans point fixe et l'image de $\gamma$ dans $\Gamma_g$
est triviale. Autrement dit,
$$\gamma \in Z({\Bbb Q}) \cap gKg^{-1} = Z({\Bbb Q}) \cap K .$$ 
Mais alors
$$gk_1 = \gamma g k_2 k' = g \gamma k_2 k'$$
et donc 
$$k_1 K' (K \cap Z({\Bbb Q})) = k_2 K' (K \cap Z({\Bbb Q}))$$
et reciproquement comme annonc\'e.

\medskip

Fixons une mesure de Haar sur $G({\Bbb A}_f )$. Alors, la famille d'applications 
\begin{eqnarray} \label{107}
\begin{array}{ccl}
H^* (Sh(G)_K ) & \rightarrow & {\Bbb C} \\
\left[ \eta \right]  & \mapsto & {\rm vol} (K/(K \cap Z({\Bbb Q}))) \cdot \int_{Sh(G)_K} \eta
\end{array}
\end{eqnarray}
d\'etermine une application 
\begin{eqnarray} \label{I}
I : H^* (Sh(G)) = \lim_{\rightarrow_K} H^* (Sh(G)_K ) \rightarrow {\Bbb C}
\end{eqnarray}
de la limite directe.

Fixons une mesure de Haar $dg$ sur $G({\Bbb A})$. Le groupe $Z({\Bbb R})$ est contenu dans $K_{\infty}$ et le quotient $K_{\infty}/Z({\Bbb R})$ est un sous-groupe
compact maximal de $G^{{\rm ad}} ({\Bbb R})^+$. De plus,
$$G({\Bbb Q}) \backslash (D \times G({\Bbb A}_f ) )/K = G({\Bbb Q}) Z({\Bbb R}) \backslash (D \times G({\Bbb A}_f ) ) / K,$$
puisque $Z({\Bbb R})$ agit trivialement sur $D$. Il correspond en particulier \`a toute forme $\eta$ sur $Sh(G)_K$ de degr\'e $pq$ (resp. bidegr\'e $(pq,pq)$) une
fonction $\widetilde{\eta}$ sur $G({\Bbb Q}) Z({\Bbb R}) \backslash G({\Bbb A})$ via les isomorphismes
\begin{eqnarray*}
\Omega^* (Sh(G)_K) & \stackrel{\sim}{\rightarrow} & \left[ \Omega^* (D) \times C^{\infty} (G({\Bbb A}_f )) \right]^{G({\Bbb Q}) \times K} \\
& \stackrel{\sim}{\rightarrow} & \left[ C^{\infty} (G({\Bbb Q}) \backslash G({\Bbb A}) ) \otimes \bigwedge \, {}^* \mathfrak{p}^* \right]^{K_{\infty} K} \\
& \stackrel{\sim}{\rightarrow} & \left[ C^{\infty} (G({\Bbb Q}) \backslash G({\Bbb A})) \right]^{K_{\infty} K},
\end{eqnarray*}
o\`u la derni\`ere application est l'\'evaluation en un vecteur orient\'e de norme un $\mathbf{1} \in \bigwedge^{pq} \mathfrak{p}$ (resp. $\in \bigwedge^{2pq} \mathfrak{p}$).
La fonction $\widetilde{\eta}$ est clairement invariante sous l'action de $Z({\Bbb R})$ et il existe une constante strictement positive $c$, ind\'ependante de $K$, telle
que
\begin{eqnarray} \label{102}
{\rm vol} (K/(K \cap Z({\Bbb Q}))) \cdot \int_{Sh(G)_K} \eta = c \cdot \int_{G({\Bbb Q}) Z({\Bbb R}) \backslash G({\Bbb A})} \widetilde{\eta} (g) dg .
\end{eqnarray}

Posons $f_{\infty} = \varphi_{\infty} (\mathfrak{1}) \in {\rm S} (V_{\infty}^p )$ fonction de Schwartz obtenue en \'evaluant en $\mathfrak{1}$ la forme de Schwartz 
$$\varphi_{\infty} \in \left[ {\rm S}(V_{\infty}^p ) \otimes \Omega^{*} (D) \right]^{G({\Bbb R})} \cong \left[ {\rm S} (V_{\infty}^p ) \otimes \bigwedge^{*} \mathfrak{p}^* \right]^{K_{\infty}},$$
\'egale \`a $\varphi^{(p)}$ en la premi\`ere place archim\'edienne et $\varphi^0_+$ en les autres.

Le th\'eor\`eme \ref{SW} implique alors le corollaire suivant.

\begin{cor} \label{CSW}
Soit $\varphi \in {\rm S} (V({\Bbb A}_f )^p )$ et $f = f_{\infty} \otimes \varphi$. Alors, il existe une constante $c' >0$, qui ne d\'epend que des choix de mesures,
telle que   
$$I([\theta (g' , \varphi )]) = c' \cdot E\left( g' , \frac12 , f \right) .$$
\end{cor}

En particulier, cette fonction est non nulle comme fonction de $g'$ (les s\'eries d'Eisenstein sont non nulles et analytiques r\'eelles en $g'$). Au vu
de la d\'ependance explicite de $[\theta (g' , \varphi )]$ en $g'$ donn\'ee par le th\'eor\`eme \ref{FD}, on en d\'eduit que 
$[\theta (g' , \varphi )]$ est non nulle pour tout $g'$ et les th\'eor\`emes \ref{T2} et \ref{T3} s'en d\'eduisent imm\'ediatement par fonctorialit\'e du cup-produit de nos
s\'eries th\^eta comme expliqu\'e en introduction \`a cette section.

\bibliography{bibliographie}

\bibliographystyle{plain}

\bigskip

\noindent
Institut de Math\'ematiques de Jussieu, \\
Unit\'e Mixte de Recherche 7586 du CNRS, \\
Universit\'e Pierre et Marie Curie, \\
4, place Jussieu 75252 Paris Cedex 05, France \\
{\it adresse electronique :} \texttt{bergeron@math.jussieu.fr}

\end{document}